\input amstex 
\input amsppt.sty    
\hsize 30pc 
\vsize 47pc 
\def\nmb#1#2{#2}         
\def\cit#1#2{\ifx#1!\cite{#2}\else#2\fi} 
\def\totoc{}             
\def\idx{}               
\def\ign#1{}             
\redefine\o{\circ} 
\define\X{\frak X} 
\define\al{\alpha} 
\define\be{\beta} 
\define\ga{\gamma} 
\define\de{\delta}

\define\et{\eta} 
\define\th{\theta} 
 
\define\ka{\kappa} 
\define\la{\lambda} 
\define\rh{\rho}

\define\ph{\varphi} 
 
\define\ps{\psi} 
\define\om{\omega} 
\define\Ga{\Gamma} 
 
\define\Th{\Theta} 
\define\La{\Lambda}

\define\Om{\Omega} 
\redefine\i{^{-1}} 
\define\x{\times} 
\define\row#1#2#3{#1_{#2},\ldots,#1_{#3}} 
\define\End{\operatorname{End}} 
\define\Fl{\operatorname{Fl}}

\define\Ad{\operatorname{Ad}} 
\define\ad{\operatorname{ad}} 
\define\pr{{\operatorname{pr}}} 
\define\Id{{\operatorname{Id}}} 
\define\tr{{\operatorname{tr}}} 
\redefine\L{{\Cal L}} 
 
\define\g{{\frak g}} 
\define\h{{\frak h}} 
\redefine\l{{\frak l}}

\def\today{\ifcase\month\or 
 January\or February\or March\or April\or May\or June\or 
 July\or August\or September\or October\or November\or December\fi 
 \space\number\day, \number\year} 
\topmatter 
\title   Poisson structures on double Lie groups 
\endtitle 
\author   
D\. Alekseevsky\\ 
J\. Grabowski\\ 
G\. Marmo\\ 
P\. W\. Michor  \endauthor 
\affil 
\leftheadtext{\smc Alekseevsky, Grabowski, Marmo, Michor} 
\rightheadtext{\smc Poisson structures} 
Erwin Schr\"odinger International Institute of Mathematical Physics,  
Wien, Austria 
\endaffil 
\address 	  
D\. V\. Alekseevsky:  
Center "Sophus Lie", gen. Antonova 2 - 99, 117279 Moscow B-279, Russia 
\endaddress 
\address J\. Grabowski: Institute of Mathematics,  
University of Warsaw,
ul. Banacha 2, PL 02-097 Warsaw, Poland; and 
Mathematical Institute, Polish Academy of Sciences,
ul. \'Sniadeckich 8, P.O. Box 137, PL 00-950 Warsaw, Poland 
\endaddress 
\email jagrab\@mimuw.edu.pl \endemail 
\address 
G. Marmo: 
Dipart\. di Scienze Fisiche - Universit\`a di Napoli, 
Mostra d'Oltremare, Pad.19, I-80125 Napoli, Italy. 
\endaddress 
\email gimarmo\@na.infn.it \endemail 
\address 
P\. W\. Michor: Institut f\"ur Mathematik, Universit\"at Wien, 
Strudlhofgasse 4, A-1090 Wien, Austria; and  
Erwin Schr\"odinger International Institute of Mathematical Physics,  
Boltzmanngasse 9, A-1090 Wien, Austria 
\endaddress 
\email peter.michor\@univie.ac.at, peter.michor\@esi.ac.at \endemail 
\date June 21, 1997 
\enddate 
\thanks  
P.W.M. was supported  
by `Fonds zur F\"orderung der wissenschaftlichen  
Forschung, Projekt P~10037~PHY'. 
\endthanks 
\abstract Lie bialgebra structures are reviewed and investigated 
in terms of the double Lie algebra, of Manin- and 
Gau{\ss}-decompositions. The standard R-matrix in a Manin 
decomposition then gives rise to several Poisson structures on the 
correponding double group, which is investigated in great detail. 
\endabstract 
\endtopmatter 
 
\document 
 
\heading Table of contents \endheading 
\noindent 1. Introduction \leaders \hbox to 
1em{\hss .\hss }\hfill {\eightrm 1}\par  
\noindent 2. Lie Bialgebras, Manin Triples, and 
Gau\ss {}-Decompositions \leaders \hbox to 
1em{\hss .\hss }\hfill {\eightrm 3}\par   
\noindent 3. Notation on Lie Groups \leaders \hbox to 
1em{\hss .\hss }\hfill {\eightrm 13}\par  
\noindent 4. Lie Poisson Groups and Double Groups \leaders \hbox to 
1em{\hss .\hss }\hfill {\eightrm 16}\par  
\noindent 5. Explicit Formulas for Poisson Structures on Double Lie 
Groups \leaders \hbox to 1em{\hss .\hss }\hfill {\eightrm 19}\par  
\noindent 6. Dressing Actions and Symplectic Leaves \leaders \hbox to 
1em{\hss .\hss }\hfill {\eightrm 27}\par  
\noindent 7. Examples \leaders \hbox to 
1em{\hss .\hss }\hfill {\eightrm 33}\par  
 
\head\totoc\nmb0{1}. Introduction \endhead 
 
In \cit!{2} we described a wide class of symplectic structures on the  
cotangent bundle $T^*G$ of a Lie group $G$ by replacing the canonical  
momenta of actions of $G$ on $T^*G$ by arbitrary ones. This method  
also worked for principal bundles and allowed us to describe the  
notion of a Yang-Mills particle which carries a `charge' given by  
spin-like variables, by means of Poisson reduction. 
 
In the latter half of this paper we consider `deformations' of $T^*G$  
in the form of so called double Lie groups equipped with the analogs  
of the symplectic structure on $T^*G$, closely related to Poisson Lie  
groups. Parts of the results may be found spread over different  
places, mainly in the unfortunately unpublished thesis of Lu  
\cit!{20}, but also to some extend in \cit!{1}, \cit!{33}, and others.  
Our presentation makes the double group the main object rather than  
Poisson Lie groups, which makes the r\^oles of $G$ and $G^*$  
manifestly symmetric and contains all the information about $G$ and  
$G^*$ and all relations between them. All this is also associated to  
the theory of symplectic groupoids as `deformed cotangent bundles' in  
general, and with mechanical systems based on Poisson symmetries as  
studied for instance in \cit!{23} and \cit!{36}. The explicit  
formulae from the second part have already found applications in  
\cit!{3}. 

The first half of this paper is devoted to the general setup:
Recall that a Poisson Lie group  is a Lie group $G$ with a Poisson
structure $\La \in \Ga(\wedge ^2TG)$ such that the multiplication map
$G\times G \to G$ is a morphism of  the Poisson manifolds.
The corresponding infinitesimal object, which determines a Poisson 
Lie group up to a covering, is that of a Lie bialgebra, defined by 
V.G.Drinfeld. It  is defined as a Lie algebra $(\g, b=[\quad,\quad])$ 
together with the structure  of a Lie algebra $(\g^*, 
b'=[\quad,\quad])$ on the 
dual space $\g^*$ such that  the bracket $b'$ defines a cocycle 
$b': \g \to \wedge^2\g$  on $\g$ with values in the $\g$-module
$\wedge^2\g$. The brackets $b,b'$ define the structure of a metrical
Lie algebra on $\l = \g \oplus \g^*$ with  Manin decomposition.
Recall that  a metrical Lie algebra is a Lie algebra together with a
non-degenerate $\ad_\g$-invariant bilinear symmetric form $g$ (the metric), 
and that a Manin  
decomposition is a decomposition of a metrical Lie algebra into 
direct sum of two isotropic subalgebras.
The metric $g$ on $\l$ is defined by the conditions that the subspaces 
$\g, \, \g^*$ are isotropic and the restriction of $g$ on 
$\g \times \g'$ is the natural pairing.
Hence, there is a natural bijection between Poisson Lie groups 
(up to a covering), bialgebras, and metrical Lie algebras with Manin
decompositions.
Remark that not every metrical Lie algebra admits a Manin decomposition
\cit!{8}.
We recall some basic constructions and facts on metrical Lie algebras
in \nmb!{2.4} -- \nmb!{2.7}.
A bivector $C \in \wedge^2\g$ on a Lie algebra $\g$ defines a 
cocycle
$$ \partial C :\g \to \wedge^2\g, X\mapsto \operatorname{ad}_XC. $$
Moreover, $C$ defines a structure of a Lie algebra on $\g^*$ if and
only if the Schouten bracket $[C,C]$ is $\operatorname{ad_\g}$ 
invariant. This condition is called the modified Yang-Baxter 
equation.

For a metrical Lie algebra $(\g,g)$ a bivector $C$ can be identified
with an endomorphism $R= C\circ g$ (the `R-matrix'). In terms of this 
endomorphism the modified Yang-Baxter equation (and other 
equations implying this) reduces to the generalized $R$-matrix 
equation (and some modifications of it), see \nmb!{2.9}.
A Manin decomposition $\g=\g_+\oplus\g_-$ of a metrical Lie algebra
$\g$ provides a solution $R= \pr_+ - \pr_-$ of the $R$-matrix equation.
More generally, we define a Gauss decomposition of a metrical Lie 
algebra $\g$ as a decomposition
$\g =\g_+ \oplus\g^0 \oplus \g_-$ of $\g$ into a sum of subalgebras such that
$\g_+, \g_-$ are isotropic and orthogonal to $\g^0$. Any solution 
$R^0$ of the $R$-matrix equation (1-mYBE) on $\g^0$, see sect 
\nmb!{2.9},
can be extended to a solution $R=\operatorname{diag}(-1,R^0,1)$ of
the same equation on $\g$. Moreover, if $R^0$ has no eigenvalues
$\pm1$, then $\g^0$ is solvable and $R^0$ is the Cayley transform
of an automorphism $A$ of $G^0$ without fixed points :
$R^0 = (A+1)(A-1)^{-1}$.
Conversly, any $R$-matrix $R$ on a metrical Lie algebra $\g$ 
defines some Gauss decomposition.

In \nmb!{2.15} we give some simple constructions of Gauss 
decompositions  of a metrical Lie algebra and its associated $R$-matrix.
Remark that the problem of describing  all bialgebra 
structures on  a given semisimple Lie algebra $\g_+$ (or the  
equivalent problem of determining all Manin decompositions
$\g =\g_+ \oplus\g_- $  of  metrical Lie algebras $\g$ with given 
$\g_+$) is solved only for a simple Lie algebra $\g_+$, 
\cit!{6}, \cit!{9}. The construction of  Weinstein of a bialgebra 
structure on a 
compact semisimple Lie algebra shows, that in general the isotropic
subalgebras $\g_+$, $\g_-$ of a Gauss decomposition of a semisimple 
Lie algebra $\g$ are not necessarily solvable. However, this is true
if the metric $g$ coincides with the Killing form of $\g$, see
\cit!{9}.

The second part of the paper is devoted to explicit description of
global versions of some objects which are studied in the first part.
The basic object is the double Lie group $G$ which corresponds to
a metrical Lie algebra $\g$ with a Manin decomposition
$\g =\g_+ \oplus\g_-$. We describe explicitly different natural Poisson and
affine Poisson structures on a double group $G$ and the dressing
action of subgroups $G_+,G_-$ associated with the isotropic 
subalgebras $\g_+, \g_-$.

\head\totoc\nmb0{2}. Lie Bialgebras, Manin Triples, and  
Gau\ss{}-Decompositions \endhead 
 
\subhead\nmb.{2.1}. Lie bialgebras and Lie Poisson groups \endsubhead 
A \idx{\it Lie bialgebra} \cit!{11} consists of a (finite dimensional)  
Lie algebra $\g$ with Lie bracket  
$b=[\quad,\quad]\in {\tsize\bigwedge}^2\g^*\otimes \g$  
and an element $b'\in{\tsize\bigwedge}^2\g\otimes \g^*$  
such that the following two properties hold: 
\roster 
\item $b'$ is a 1-cocycle $\g\to{\tsize\bigwedge}^2\g$:  
     $\partial_bb'=0$ where   
     $(\partial_bb')(X,Y)=-b'([X,Y])+\ad_X(b'(Y))-\ad_Y(b'(X))$. 
     To put this into perspective, note that this is equivalent to  
     the fact that $X\mapsto (X,b'(X))$ is a homomorphism of Lie  
     algebras from $\g$ into the semidirect product  
     $\g\ltimes {\tsize\bigwedge}^2\g$ with the Lie bracket  
     $[(X,U),(Y,V)] = ([X,Y],\ad_XV-\ad_YU)$.  
\item $b'$ is a Lie bracket on $\g^*$. 
\endroster 
In \cit!{17} a graded Lie bracket on ${\tsize\bigwedge}(\g\x \g^*)$  
is constructed  
which recognizes Lie bialgebras, their representations,  
and gives the associated notion of Chevalley cohomology.

\subhead\nmb.{2.2}. Exact Lie bialgebras and Yang-Baxter  
equations \endsubhead 
A Lie bialgebra $(\g,b,b')$ is called \idx{\it exact} if the  
1-cocycle $b'$ is a coboundary:  
$b'=\partial_bC$ for $C\in{\tsize\bigwedge}^2\g$,   
i.e\. $b'(X)=\ad_XC$. 
A bivector $C\in{\tsize\bigwedge}^2\g$ defines a Lie bialgebra structure  
$b'=\partial_bC$ on $\g$ if and only if the Schouten bracket (see  
\nmb!{3.4}) is  
$\ad(\g)-invariant$: 
$$ 
[C,C]\in ({\tsize\bigwedge}^3\g)^\g. \tag{mYBE} 
$$ 
This condition is called the \idx{\it modified Yang-Baxter Equation}. 
In particular any Poisson bivector $C\in {\tsize\bigwedge}^2\g$ satisfying 
$$ 
[C,C] = 0 \tag{YBE} 
$$ 
defines a bialgebra structure $b'=\partial_bC$ in $\g$. This equation  
is called the \idx{\it Yang-Baxter Equation}. 
 
If $\g$ is semisimple then by the Whitehead lemma  
$H^1(\g,{\tsize\bigwedge}^2\g)=0$, so any cocycle $b'$ is a coboundary, and the  
classification of all bialgebra structures on $\g$ reduces to the  
description of all bivectors $C\in{\tsize\bigwedge}^2\g$ which satisfy  
\thetag{mYBE}. If moreover the Lie algebra $\g$ is simple then the  
space $({\tsize\bigwedge}^3\g)^\g$ is 1-dimensional, generated by the 3-vector  
$B^g\in{\tsize\bigwedge}^3\g$ given by 
$B^g(\al,\be,\ga):= g([g\i\al,g\i\be],g\i\ga)$, where $g$ denotes the  
Cartan-Killing form. So for simple $\g$ the modified Yang-Baxter Equation  
\thetag{mYBE} can be written, using the Schouten bracket, as 
$$ 
[C,C] = cB^g 
$$ 
All solutions of this equation for $c\ne0$ for complex simple $\g$  
were described by \cit!{6}, \cit!{9}. 
 
\subhead\nmb.{2.3}. Manin decompositions \endsubhead 
Let $(\g,b)$ be a Lie algebra and let $b'$ be a Lie bracket on the  
dual space $\g^*$. Let us define a skew symmetric bracket  
$[\quad,\quad]$ on the vector space $\l:=\g\oplus\g^*$ by  
$$ 
[(X,\al),(Y,\be)] := \Bigl(b(X,Y)+\ad_{b'}^*(\al)Y-\ad_{b'}^*(\be)X, 
          b'(\al,\be)+\ad_{b}^*(X)\be-\ad_{b}^*(Y)\al\Bigr), 
$$ 
where $\ad_b(X)Y=b(X,Y)$, $\ad_b^*(X)=\ad_b(-X)^*\in \End(\g^*)$, and  
similarly for $b'$. The adjoint operator $\ad(X,\al)\in \End(\l)$  
is skew symmetric with respect to the natural pseudo Euclidean inner  
product $g$ on $\l$ which is given by  
$g((X,\al),(Y,\be))=\langle \al,Y\rangle + \langle \be, X\rangle$,  
and the skew symmetric bracket is uniquely determined by this  
property. 
The skew symmetric bracket $[\quad,\quad]$ on $\l$ satisfies the  
Jacobi identity if and only if $b':\g\to{\tsize\bigwedge}^2\g$ is a 1-cocycle  
with respect to $b$: $\partial_bb'=0$; or equivalently if and only if  
$b:\g^*\to {\tsize\bigwedge}^2\g^*$ is a 1-cocycle with respect to $b'$:  
$\partial_{b'}b=0$.  
 
Following Astrakhantsev \cit!{4} we will call 
\idx{\it metrical Lie algebra} a Lie algebra $\l$ together with  
an $\ad$-invariant inner product $g$: $g([X,Y],Z)=g(X,[Y,Z])$. 
 
A decomposition of a metrical Lie algebra $(\l,g)$ as a direct sum  
$\l=\g_+\oplus\g_-$ of two $g$-isotropic Lie subalgebras $\g_+$ and  
$\g_-$ is called a \idx{\it Manin decomposition}.  
 
A triple of Lie algebras $(\g,\g_+,\g_-)$ together with a duality  
pairing between $\g_+$ and $\g_-$ is called a \idx{\it Manin triple}  
if $\g=\g_+\oplus\g_-$, $\g_+$ and $\g_-$ are Lie algebras of $\g$,  
and the duality pairing induces an $ad$-invariant inner product on $\g$  
for which $\g_+$ and $\g_-$ are isotropic. 
 
\proclaim{Theorem} \cit!{10} There exist a natural bijective  
correspondence between Lie bialgebras $(\g,b,b')$ and metrical Lie  
algebras $(\frak l,g)$ with Manin decomposition  
$\frak l= \g\oplus\g^*$.  
\endproclaim 
 
The Lie algebra $\frak l=\g\oplus\g^*$ associated to the Lie  
bialgebra $(\g,b,b')$ is called the \idx{\it Manin double}.  
 
\subhead\nmb.{2.4}. Examples of metrical Lie algebras \endsubhead 
Any commutative Lie algebra has the structure of a metrical Lie  
algebra, with respect to any inner product. Any semisimple Lie algebra  
is metrical, the metric is given by the Cartan-Killing form. 
 
Let $\g$ be a Lie algebra. Let us denote by $T^*\g=\g\ltimes \g^*$ the  
semidirect product of the Lie algebra $\g$ with the abelian ideal  
$\g^*$, where $\g$ acts on $\g^*$ by the the coadjoint action. This  
is the Lie algebra of the cotangent group $T^*G$ of a Lie group $G$  
with Lie algebra $\g$. The natural pairing between $\g$ and the dual  
$\g^*$ defines an $\ad$-invariant inner product $g$ on $T^*\g$ for  
which the subalgebras $\g$ and $\g^*$ are isotropic, by definition of  
the coadjoint action. Hence $T^*\g=\g\oplus \g^*$ is a Manin  
decomposition of the metrical Lie algebra $T^*\g$. 
It describes the Lie bialgebra structure $b'=0$ on $\g$. 
 
\subhead\nmb.{2.5} \endsubhead 
We will now describe the \idx{\it double extension} of a metrical Lie  
algebra according to Kac \cit!{14},~2.10, Medina-Revoy \cit!{25}:  
Let $(\g,g)$ be a metrical Lie algebra and let $\frak d$ be a Lie  
algebra together with a representation  
$\rh:\frak d\to \operatorname{Der}_{\text{skew}}(\g,g)$ by skew  
symmetric derivations on $\g$. We then put   
$$\align 
\g_{\frak d} &:= \frak d \oplus \g \oplus \frak d^*,\\ 
[D_1+X_1+\al_1,D_2+X_2+\al_2] &= [D_1,D_2]_{\frak d} + \\  
&\quad + [X_1,X_2]_{\frak g} + \rh(D_1)(X_2) - \rh(D_2)(X_1) + \\ 
&\quad + c(X_1,X_2) + \ad_{\frak d}^*(D_1)(\al_2) -  
\ad_{\frak d}^*(D_2)(\al_1),\\ 
g_{\g_{\frak d}}(D_1+X_1+\al_1,D_2+X_2+\al_2) &= g(X_1,X_2) +  
\langle \al_1,D_2 \rangle + \langle \al_2, D_1 \rangle, 
\endalign$$ 
where the central cocycle  
$c:\g\x\g\to \frak d^*$ is given by  
$\langle D,c(X,Y) \rangle = g(\rh(D)(X),Y)$ for $D\in\frak d$.  
Then $\g_{\frak d}$ is again a metrical Lie  
algebra. Note that the metrical Lie subalgebra  
$\frak d\oplus \frak d^*$ is isomorphic to the cotangent Lie algebra  
$T^*\frak d$ and that we may view $\g_{\frak d}$  
as the semidirect product $\g_{\frak d} = \frak d\ltimes \h$, where  
$\h$ is the central extension 
$$ 
0\to \frak d^* \to \h \to \g \to 0 
$$ 
described by the cocycle $c$ and where $\frak d$ acts on $\h$ by  
$(\rh,\ad_{\frak d}^*)$. 
 
The orthogonal direct sum of two metrical Lie algebras is again a  
metrical Lie algebra. In particular the orthogonal direct sum of a  
metrical Lie algebra $\g$ with a 1-dimensional abelian metrical Lie  
algebra is called the \idx{\it trivial extension} of $\g$.  
 
\proclaim{Theorem} Kac \cit!{14},~2.11, Revoy-Medina, \cit!{25}. 
Any solvable metrical Lie algebra can be obtained from a commutative  
metrical Lie algebra by an appropriate sequence of double extensions  
and trivial extensions. 
\endproclaim 
 
\subhead\nmb.{2.6} \endsubhead 
The following result gives an analogon of the Levi-Maltsev  
decomposition for a metrical Lie algebra. 
 
\proclaim{Theorem} Astrachantsev \cit!{4}.  
Any metrical Lie algebra  
$\g$ is an orthogonal direct sum  
$$ 
\g = \h\oplus \frak r = \frak s_n\oplus T^*\frak s_i \oplus \frak r 
$$ 
consisting of a subalgebra $\h$ with commutative radical and a  
solvable ideal $\frak r$. Moreover, $\h$ is an orthogonal direct sum  
of a maximal $g$-non-degenerate semisimple Lie subalgebra $\frak s_n$  
and the cotangent algebra $T^*\frak s_i$ of a maximal $g$-isotropic  
semisimple Lie subalgebra $\frak s_i$ of $\g$.  
\endproclaim 
 
\subhead\nmb.{2.7}.  Metrical extensions \endsubhead 
Bordemann \cit!{8} gave the following construction of a metrical Lie  
algebras. 
 
Let $\frak a$ be a Lie algebra and let  
$w: \frak a\wedge \frak a\to \frak a^*$ be a 2-cocyle with 
values in the $\frak a$-module $\frak a^*$.  
Then the Lie algebra extension  
$$0\to \frak a^*\to \g_w\to \frak a\to 0$$  
described by $w$, i\.e\. the Lie algebra  
$\g_w:=\frak a\oplus\frak a^*$ with bracket 
$$[(a,\al),(b,\be)]_{\g_w}:= ([a,b]_\g,w(a,b)+\ad^*(a)\be-\ad^*(b)\al)$$ 
is a metrical Lie algebra with metric 
$$g((a,\al),(b,\be)):= \langle \al,b \rangle+\langle\be,a \rangle$$ 
if and only if $w$ has the following property 
$$ 
\langle a,w(b,c) \rangle = \langle b,w(c,a) \rangle \quad 
\text{ for }a,b,c\in\frak a. 
$$  
If $w=0$ then this is exactly the metrical Lie algebra $T^*\frak a$ ---  
thus Bordemann called this construction the $T^*$-extension. 
 
\proclaim{Theorem} \cit!{8} Any $2n$-dimensional complex solvable  
metrical Lie algebra $\g$ is a metrical extension of some  
$n$-dimensional Lie algebra $\frak a$. Moreover any isotropic ideal  
of $\g$ is contained in an $n$-dimensional isotropic commutative ideal  
of $\g$.  
\endproclaim 
   
\subhead\nmb.{2.8}. The Yang-Baxter equations on metrical Lie algebras 
\endsubhead 
In the case of a metrical Lie algebra $(\g,g)$ we can pull down one  
index of bivector $C\in {\tsize\bigwedge}^2 g$ and we can  
reformulate the (modified) Yang-Baxter equation in terms of the operator  
$R=C\o g:\g\to\g^*\to \g$. 
 
First let $(\g,b=[\quad,\quad])$ be a Lie algebra. For any  
$R\in \End(\g)$ we   
define two elements $b_R, B_R\in \g\otimes {\tsize\bigwedge}^2\g^*$ by  
$$\gather 
b_R(X,Y) = [X,Y]_R := [RX,Y]+[X,RY], \\ 
B_R(X,Y) := [RX,RY] -R[X,Y]_R = [RX,RY] - R[RX,Y] - R[X,RY].  
\endgather$$ 
Note that $B_R$ is related to the Fr\"olicher-Nijenhuis-like bracket  
$[R,R]$ by  
$$\align 
\frac12[R,R](X,Y)&=[RX,RY]-R([RX,Y]+[X,RY])+R^2[X,Y] =\\  
&=B_R(X,Y)+R^2[X,Y]. 
\endalign$$ 
 
\proclaim{Proposition} 
Let $(\g,g)$ be a metrical Lie algebra, let $C\in{\tsize\bigwedge}^2\g$ and let  
$R=C\o g:\g\to\g^*\to \g$ be the corresponding operator.  
Then we have: 
\roster 
\item Via the isomorphism $g\i:\g^*\to \g$  
     the bracket $b'=\partial_bC\in \g^*\otimes {\tsize\bigwedge}^2\g$ on  
     $\g^*$ corresponds to the bracket $b_R$ on $\g$: 
$$ 
g\i(b'(\al,\be)) = b_R(g\i\al,g\i\be) = [g\i\al,g\i\be]_R, 
\quad\text{ for }\al,\be\in\g^*. 
$$ 
\item Under the embedding  
     ${\tsize\bigwedge}^3\g\to \g\otimes{\tsize\bigwedge}^2\g^*$ induced   
     by $g$, the Schouten bracket $[C,C]\in {\tsize\bigwedge}^3\g$  
     corresponds to the element  $2B_R\in\g\otimes {\tsize\bigwedge}^2\g^*$. 
\endroster 
\endproclaim 
 
\demo{Proof}  
Let $X,Y,Z\in \g$ and $\al=gX,\be=gY\in\g^*$.  
Note that $g(RX,Y)=g(X,-RY)$. Then  
$$\align 
\langle Z, b'(\al,\be)\rangle &= \langle \ad_ZC,\al\wedge \be\rangle   
     = \langle C,(\ad_Z)^*\al\wedge \be + \al\wedge (\ad_Z)^*\be \rangle\\ 
&= \langle C,(\ad_Z)^*gX\wedge gY + gX\wedge (\ad_Z)^*gY \rangle\\ 
&= \langle C,-g\ad_ZX\wedge gY - gX\wedge g\ad_ZY \rangle\\ 
&= -\langle Cg\ad_ZX, gY \rangle - \langle CgX, g\ad_ZY \rangle\\ 
&= -g(R[Z,X],Y) - g(RX,[Z,Y]) =  g(Z,[X,RY]) + g([RX,Y],Z) \\ 
&= \langle Z ,g[X,Y]_R\rangle. 
\endalign$$ 
 
For proving the second assertion we may assume  
without loss that $C\in{\tsize\bigwedge}^2\g$ is  
decomposable, $C=X\wedge Y$, since both sides are quadratic.  
Then we have:  
$$\align 
R(Z)&=(C\o g)(Z)=((X\wedge Y)\o g)(Z) = g(Y,Z)X-g(X,Z)Y\\ 
B_R(U&,V) = [RU,RV] - R[RU,V] - R[U,RV] \\ 
&= [g(Y,U)X-g(X,U)Y,g(Y,V)X-g(X,V)Y] \\ 
&\quad - g(Y,[g(Y,U)X-g(X,U)Y,V])X+g(X,[g(Y,U)X-g(X,U)Y,V])Y\\ 
&\quad- g(Y,[U,g(Y,V)X-g(X,V)Y])X+g(X,[U,g(Y,V)X-g(X,V)Y])Y\\ 
&= -g(Y,U)g(X,V)[X,Y] -g(X,U)g(Y,V)[Y,X] \\ 
&\quad - g(Y,U)g([Y,X],V)X - g(X,U)g([X,Y],V)Y\\ 
&\quad + g(Y,V)g([Y,X],U)X + g(X,V)g([X,Y],U)Y\\ 
\endalign$$ 
On the other hand we have for the Schouten bracket  
$$\align 
[C,C] &= [X\wedge Y,X\wedge Y]  
     = 2[X,Y]\wedge X\wedge Y \\ 
\frac12\langle [C,C]&, \al\wedge gU\wedge gV \rangle  
     = \langle [X,Y]\wedge X\wedge Y, \al\wedge gU\wedge gV \rangle \\ 
&=\det\pmatrix \langle [X,Y],\al \rangle & \langle X,\al \rangle &  
                              \langle Y,\al \rangle\\ 
               g([X,Y],U) & g(X,U) & g(Y,U)\\ 
               g([X,Y],V) & g(X,V) & g(Y,V)\endpmatrix \\ 
&=\langle B_R(U,V),\al \rangle, 
\endalign$$ 
from the computation above. 
\qed\enddemo 
 
\subhead Remarks \endsubhead 
We may extend $R\mapsto B_R$ to a bracket in  
${\tsize\bigwedge}\g^*\otimes\g$ as follows.  
On decomposable tensors this bracket is given by  
$$ 
[\ph\otimes X, \ps\otimes Y] = \ph \wedge \ps \otimes [X,Y]  
+ \ph \wedge \ad^*_X\ps\otimes Y - \ad^*_Y\ph \wedge \ps \otimes X, 
$$ 
and it defines a $\Bbb Z$-graded Lie bracket on  
${\tsize\bigwedge}^*\g^*\otimes \g$. If $\g$ acts by derivations on a graded  
commutative algebra $A=\bigotimes_{i=0}^\infty A_i$, the same formulae  
define a graded Lie bracket on $A\otimes \g$. 
 
Moreover we have $B_R=\frac12[R,R]^B$, and by the graded Jacobi  
identity we get the analogon of the Bianchi identity 
$[R,B_R]^B=0.$ 
 
The invariant inner product $g:\g\to\g^*$ induces an embedding 
$$\tsize{\tsize\bigwedge}^{*+1}\g \to \tsize{\tsize\bigwedge}^*\g^*\otimes\g$$ 
which is a homomorphism from the Schouten bracket to the graded Lie  
bracket $[\quad,\quad]^B$. This follows from the polarization of  
\therosteritem2 in the proposition above (note that the brackets in  
degree 1 are symmetric), since $\g$ and ${\tsize\bigwedge}^2\g$ generate the  
whole Schouten algebra.  
 
On a manifold one may also consider the bracket $[\quad,\quad]^B$ but  
it maps tensor fields to differential operators.  
 
There is a homomorphism of graded Lie algebras  
$$\align 
({\tsize\bigwedge}^*\g^*\otimes \g , [\quad,\quad]^B)  
     &\to (\Om^*(\g,\g),[\quad,\quad]^{FN}),\\ 
\al_1\wedge \dots\wedge \al_p \otimes X  
     &\mapsto d\al_1\wedge \dots\wedge d\al_p \otimes \ad_{\g}(X), 
\endalign$$ 
where $\Om(\g,\g)\cong \Om(\g;T\g)$ is the graded Lie algebra of all  
tangent space valued differential forms on $\g$ with the  
Fr\"olicher-Nijenhuis bracket. The kernel of this homomorphism  
consists of ${\tsize\bigwedge}^*\g^*\otimes Z(\g)$ where $Z(\g)$ is  
the center of $\g$. All this follows from the well known formula  
for the Fr\"olicher-Nijenhuis bracket (see e.g\. \cit!{16},~8.7) 
$$\align 
[\ph\otimes \xi, \ps\otimes \et]  
&= \ph \wedge \ps \otimes [\xi,\et] + \ph \wedge 
\L_\xi\ps\otimes \et - \L_\et\ph \wedge \ps \otimes \xi \\ 
&\quad + (-1)^{\deg\ph}\left(d\ph \wedge i_\xi\ps \otimes \et + i_\et\ph 
\wedge  d\ps \otimes \xi \right), 
\endalign$$ 
where $\ph,\ps\in\Om(\g)$ are differential forms and where  
$\xi,\et\in\X(\g)$ are vector fields. 
 
 
\proclaim{\nmb.{2.9}. Corollary} (see \cit!{32}) 
For $C\in {\tsize\bigwedge}^2\g$ and $R=C\o g:\g\to\g$ the following  
conditions are equivalent. 
\roster 
\item $b'=\partial C$ is a Lie bracket in $\g^*$, hence $(\g,b,b')$  
       is a Lie bialgebra. 
\item $b_R$ is a Lie bracket in $\g$. 
\item The Schouten bracket $[C,C]\in{\tsize\bigwedge}^3\g$ is  
       $\ad_\g$-invariant.  
\item $B_R\in (\g\otimes{\tsize\bigwedge}^2\g^*)$ is $\g$-invariant.  
\item For all $X,Y,Z\in\g$ we have 
$$ 
[X,B_R(Y,Z)]+[Y,B_R(Z,X)]+[Z,B_R(X,Y)]=0.  
$$ 
\endroster 
\endproclaim 
 
\demo{Proof} 
It remains to show that \therosteritem4 is equivalent to  
\therosteritem5. This follows from the identity 
$$ 
g((\ad(U)B_R)(Y,Z),X) = -g([X,B_R(Y,Z)]+[Y,B_R(Z,X)]+[Z,B_R(X,Y)],U) 
$$ 
which holds for all $X,Y,Z,U\in\g$. 
\qed\enddemo 
 
The following simpler equations obviously imply equation  
\therosteritem5: 
$$\gather 
B_R+I\o b=0,\quad\text{ or }\quad B_R(X,Y)+I[X,Y]=0, \tag{I-mYBE}\\ 
B_R+cb=0,\quad\text{ or }\quad B_R(X,Y)+c[X,Y]=0, \tag{c-mYBE}\\ 
[C,C]=0\quad\text{ or }\quad B_R=0\tag{YBE} \\ 
\endgather$$ 
where $I\in \End(\g)^\g$ is an $\ad_\g$-invariant operator on $\g$, and  
where $c$ is a constant in $\Bbb K$.  
If $\Bbb K=\Bbb C$ (or $\Bbb K=\Bbb R$) without loss we may  
assume that $c=1$ (or $c=\pm1$).  
 
In \cit!{9}, 3.2 it was shown that any structure of a bialgebra on a  
semisimple Lie algebra comes from a solution of \thetag{I-mYBE} for  
some $I\in\End(\g)^\g$; and for a simple Lie algebra from a solution  
of \thetag{c-mYBE}.  
 
It is also interesting to construct non skew symmetric solutions of  
all this equations. Some class of solutions on a simple complex Lie  
algebra was constructed by \cit!{31}. 
 
Note that for an $\ad_\g$-invariant operator $I\in\End(\g)^\g$ we have 
$B_I= I^2\o b$ since $I[X,Y]=[IX,Y]=[X,IY]$. So any skew symmetric  
$\ad$-invariant operator $I$ gives a solution of the \thetag{mYBE}. 
Nonconstant operators of this kind exist on semisimple Lie algebras  
$\g$ if and only if $\g$ has isomorphic simple summands: For example,  
if $\g=l\g_1 = \g_1\oplus\dots\oplus\g_1=\g_1\otimes \Bbb K^l$ then  
$\End(\g)^\g=1\otimes \End(\Bbb K^l)$, and any skew symmetric  
matrix $A\in\End(\Bbb K^l)$ gives a solution $I=1\otimes A$ of  
\thetag{mYBE}. 
 
 
To distinguish equations for $C\in{\tsize\bigwedge}^2\g$ and for $R=C\o g$ the  
equation \thetag{1-mYBE} for  
$R$ will be called the \idx{\it $R$-matrix equation}, and solutions  
will be called \idx{\it $R$-matrices}. 
 
\subhead\nmb.{2.10} \endsubhead
Let $(\g,b,g)$ be a metrical Lie-algebra and let  
$R\in\operatorname{End}(\g)$ be a skew symmetric endomorphism. 
 
\proclaim{Lemma} \cit!{32}, \cit!{9} 
The following conditions are equivalent. 
\roster 
\item The endomorphism $R$ satisfies the $R$-matrix equation  
       $B_R+b=0$. 
\item The endomorphisms $R_\pm:=R\pm1$ satisfy 
$$ 
R_+[R_-X,R_-Y] = R_-[R_+X,R_+Y]\quad\text{ for }X,Y\in\g. 
$$ 
\item For all $\la,\mu\in \Bbb C$ and $X,Y\in\g$ we have 
$$\align
(\la+\mu)R[X,Y] = &(1+\la\mu)[X,Y] + 
+ [(R-\la)X,(R-\mu)Y] \\
&-  (R-\la)[X,(R-\mu)Y] - (R-\mu)[(R-\la)X,Y]. 
\endalign$$ 
\item The bracket $ b_R(X,Y) = [X,Y]_R = [RX,Y] + [X,RY]$ is a Lie  
       bracket and moreover both $R_\pm:(\g,b_R)\to (\g,b)$ are  
       homomorphisms.  
\endroster 
\endproclaim 
 
\subhead\nmb.{2.11} \endsubhead 
For an endomorphism $R:\g\to\g$ and $\la\in\Bbb C$  
the space 
$$ 
\g_\la=\ker(R-\la)^N\quad\text{ for large }N 
$$ 
is called \idx{\it weight space} if it is not 0, and $\la$ is called  
weight of $R$. We have the following decomposition of $\g$ into a  
direct sum of all weight spaces 
$$ 
\g = \bigoplus_{\la\in W}\g_\la, 
$$ 
where $W$ is the set of all weights. 
 
For $\la,\mu\in \Bbb C$ with $\la+\mu\ne0$ we put  
$$ 
\la\circ\mu:= \frac{1+\la\mu}{\la+\mu}.\tag1 
$$ 
Note that $(\pm1)\circ\mu=\pm1$. 
 
\proclaim{Lemma} \cit!{9} Let $R$ be an $R$-matrix on  
a metrical Lie algebra $(\g,g)$. Then we have: 
\roster 
\item For weights $\la$, $\mu$ with $\la+\mu\ne0$ we have 
$$ 
[\g_\la,\g_\mu]\subseteq\g_{\la\circ\mu}\quad
\text{ and }\quad g(\g_\la,\g_\mu)=0. 
$$ 
\item For $\la\ne\pm1$ we have $[\g_\la,\g_{-\la}] =0$. 
\item The spaces $\g_{\pm1}$ are Lie subalgebras of $\g$, and  
     $[\g_\la,\g_{\pm1}]\subseteq\g_{\pm1}$ for $\la\ne\pm1$. 
\endroster 
\endproclaim

\subhead\nmb.{2.12}. $R$-matrices and associated Gau\ss{}  
decompositions \endsubhead  
We will discuss the relations between $R$-matrices on a metrical Lie  
algebra and its Gau\ss{} decompositions. 
 
\subhead Definition \endsubhead 
A \idx{\it (generalized) Gau\ss{} decomposition} of a metrical Lie  
algebra $(\g,g)$ is a decomposition of $\g$  
$$ 
\g=\g_+\oplus\g^0\oplus\g_- 
$$ 
into a sum of subalgebras, where the inner product $g$ is non  
degenerate on $\g_0$, and where $\g_+$ and $\g_-$ are isotropic  
subalgebras which are orthogonal to $\g_0$.  
 
Note that a Manin decomposition is the special case of a Gau\ss{}  
decomposition with $\g^0=0$.  
 
\proclaim{Proposition} 
An $R$-matrix $R$ on a metrical Lie-algebra $(\g,g)$ defines a  
Gau\ss{} decomposition  
$$ 
\g = \g_-\oplus \g^0 \oplus \g_+, 
$$ 
where $\g_\pm$ are the weight spaces $\g_{\pm1}$ of $R$, and where  
$$ 
\g^0=\bigoplus_{\la\ne\pm1}\g_\la 
$$  
is a solvable Lie subalgebra which admits an $g$-orthogonal  
automorphism $A=((R+1)|\g^0)\o ((R-1)|\g^0)\i$ without fixed point (so  
$AX=X$ implies $X=0$). 
 
Conversely, let $\g = \g_-\oplus \g^0 \oplus \g_+$ be a Gau\ss{}  
decomposition of a metrical Lie algebra $(\g,g)$, where $g^0$ admits  
an orthogonal automorphism $A$ without fixed points.  
Put $R_0=(A+1)\o(A-1)\i$. Then  
$$ 
R=\operatorname{diag}(-1,R_0,1):\g\to\g 
$$ 
is an $R$-matrix.  
 
More generally, any $R$-matrix $R'$ on $\g$ which induces this Gauss{}  
decomposition has the form  
$$ 
R'=\operatorname{diag}(-1+N_-,R_0,1+N_+):\g\to\g,  
$$ 
where $N_\pm:\g_\pm\to\g_\pm$ are suitable nilpotent endomorphisms. 
\endproclaim 
 
Remark that, 
in fact, the $R$-matrix equation specifies the form of $N_\pm$. For  
example, denote by $\g^i_\pm = \ker(N_\pm)^i\subseteq \g_\pm$.  
Then 
$$ 
\g_\pm =\g_\pm^k \supset \g_\pm^{k-1}\supset \dots \supset \g_\pm^{1}  
\supset \g_\pm^{0} = 0 
$$ 
is a chain of ideals: $[\g_\pm^{i},\g_\pm]\subset \g_\pm^{i}$.  
 
\demo{Proof} 
The first statement follows immediately from lemma \nmb!{2.11}. 
The operators $R_\pm|\g^0$ are invertible.   
Note that by putting $X=(R-1)\i u$ and $Y=(R-1)\i v$ for $u,v\in\g^0$  
the equation in lemma \nmb!{2.10}.\therosteritem2 becomes 
$(R+1)(R-1)\i[u,v]=[(R+1)(R-1)\i u,(R+1)(R-1)\i v]$. 
This shows that $A=(R+1)(R-1)\i$ is an automorphism of $\g^0$.  
It has no fixed point. It is easily seen that $A$ is orthogonal if  
and only if $R|g^0$ is skew symmetric.  
 
We now use the fact that a Lie algebra which admits an automorphism  
without fixed point is solvable, see \cit!{37}.  
 
For the converse, since all arguments above were equivalencies, we see  
that $R_0=(A+1)(A-1)\i$ is a (skew symmetric) $R$-matrix on $\g^0$.  
Using lemma \nmb!{2.10}.\therosteritem2 again it follows by  
checking cases $X,Y\in\g_-,\g_+,\g_0$  
that $R=\operatorname{diag}(-1,R_0,1)$ is an $R$-matrix. 
 
The last statement is obvious. 
\qed\enddemo 
 
\proclaim{\nmb.{2.13}. Corollary} 
Any semisimple $R$-matrix $R$ on a metrical Lie algebra $(\g,g)$ can  
be written as  
$$ 
R=\operatorname{diag}(-1,R_0,1) 
$$ 
with respect to an appropriate Gau\ss{} decomposition 
$\g = \g_-\oplus \g^0 \oplus \g_+$,  
where $R_0= (A+1)(A-1)\i$ for a semisimple orthogonal automorphism $A$ of  
$\g_0$ without fixed point. \qed 
\endproclaim 
 
\proclaim{\nmb.{2.14}. Corollary} 
Any $R$-matrix $R$ on a metrical Lie algebra $(\g,g)$ without  
eigenvalues $\pm1$ is of the form 
$$ 
R=(A+1)\o(A-1)\i, 
$$ 
where $A$ is an orthogonal automorphism of $\g$ without fixed point.  
\qed 
\endproclaim 
 
Note that non-orthogonal automorphisms $A$ give non-skew symmetric  
solutions   
$R=(A+1)\o (A-1)\i$ 
of the $R$-matrix equation.  
 
\subhead\nmb.{2.15}. Construction of $R$-matrices via Gauss{}  
decompositions \endsubhead 
Let $(\g,g)$ be a metrical Lie algebra.  
Choose a skew-symmetric  
derivation $D$ of $\g$ (for example an inner derivation $\ad(X_0)$  
for $X_0\in\g$).  
It defines a decomposition  
$$\gather 
\g=\g_-\oplus\g^0\oplus g_+,\quad 
     \text{ where }\quad \g^0=\g_0\quad\text{ and }\\ 
\g_+=\negthickspace\bigoplus\Sb\Re(\la)>0\text{ or }\\\Re(\la)=0, 
                    \Im(\la)>0\endSb\negthickspace\g_\la, \qquad 
\g_-=\negthickspace\bigoplus\Sb\Re(\la)>0\text{ or }\\\Re(\la)=0, 
                    \Im(\la)>0\endSb\negthickspace\g_{-\la}. 
\endgather$$ 
 
\proclaim{Lemma} 
For an complex Lie algebra this decomposition associated to a skew  
symmetric derivation $D$ is a Gau\ss{} decomposition. 
\endproclaim 
 
\demo{Proof} 
$g((D-\mu)^lX,Y)=g(X,(-D-\mu)^lY)$. 
\qed\enddemo  
 
We can iterate this construction if there exists non-nilpotent skew  
symmetric derivations of $\g_0$, in particular if $\g_0$ is not  
nilpotent. Hence we have: 
 
\proclaim{\nmb.{2.16}. Corollary} 
Let $D$ be a skew symmetric derivation on $(\g,g)$. 
 
The decomposition associated to $D$ is trivial, $\g=\g^0$, if and  
only if $D$ is nilpotent.  
 
If 0 is not an eigenvalue of $D$ then the associated decomposition is  
a Manin-decomposition 
$$ 
\g=\g_+\oplus\g_-. 
$$ 
\endproclaim  
 
 
\subhead\nmb.{2.18}. Remark \endsubhead 
In the special case when the subalgebra $\g^0$ of a Gauss{}  
decomposition is commutative, then for any skew symmetric  
endomorphism $R_0:\g^0\to\g^0$ the operator 
$$ 
R=\operatorname{diag}(-1,R_0,1) 
$$ 
is an $R$-matrix. 
It is known, \cit!{12}, or \cit!{24},~9.3.10, 
that the connected component of the stabilizer 
of a regular point in the coadjoint representation of any connected 
Lie group is   
commutative. For a metrical Lie algebra the adjoint representation is  
isomorphic to the coadjoint one. Hence the Gau\ss{} decomposition  
associated to an inner derivation $\ad(X)$ of a regular semisimple element  
$X\in\g$ has $\g^0$ commutative.  
 
\subhead\nmb.{2.19}. Construction of $R$-matrices without  
eigenvalues $\pm1$ \endsubhead 
Let $\g$ be a (nilpotent) Lie algebra which admits a derivation with  
positive eigenvalues. For example, let $\g=\bigoplus_{i>0}\g_i$ be a  
positively graded Lie algebra and let $D|\g_i=i\Id$. 
Denote by $T^*\g=\g\ltimes \g^*$ the semidirect sum of $\g$ and the  
commutative ideal $\g^*$ with the coadjoint action on $\g^*$.  
The natural pairing $\g\x\g^*\to\Bbb C$ defines an  
$\ad_{T^*\g}$-invariant metric $g$ on $\g$. The derivation $D$ can  
naturally be extended to a $g$-skew symmetric derivation $D$ on  
$T^*\g$ without eigenvalue 0.  
Then $A_t:=\exp(tD)$ is a $g$-orthogonal automorphism of $(T^*\g,g)$  
without fixed point. Hence  
$$ 
R=(A_t+1)(A_t-1)\i 
$$ 
is an $R$-matrix without eignevalues $\pm1$. 
 
\proclaim{\nmb.{2.20}. Proposition} \cit!{32} Let  
$\g=\g_+\oplus\g_-$ be a Manin decomposition  
of a metrical Lie algebra $\g$, and let  
$\pr_{\pm}:\g\to\g_{\pm}$ be the corresponding projections. 
Then $R=\pr_+-\pr_-$ is a solution of \thetag{1-mYBE} $B_R+b=0$. 
\endproclaim 
 
\proclaim{\nmb.{2.21}. Proposition} Let  
$\g=\g_+\oplus\g_0\oplus\g_-$ be a Gau\ss{}-decomposition of a  
metrical Lie algebra $\g$, and let $\pr_{\pm}:\g\to\g_\pm$ be the  
orthogonal projections. Then any solution $R_0$ of the  
\thetag{1-mYBE} on $\g_0$ has an extension  
$R=c(\pr_+\oplus R_0\oplus(1-c)\pr_-)$ to a solution of the  
\thetag{1-mYBE} on $\g$.  
\endproclaim 
 
This gives us an inductive procedure for the construction of  
solutions of the \thetag{mYBE}. 
 
\proclaim{\nmb.{2.22}. Theorem} 
Let $(\g,g)$ be a metrical Lie algebra and let $R:\g\to\g$ be a  
solution of \nmb!{2.9}, \thetag{1-mYBE}. 
Then the following Manin decompositions are isomorphic: 
\roster 
\item The Manin double $\g\oplus\g^*$ associated to the bialgebra  
     structure $b'=\partial_b(R\o g\i)$ from \nmb!{2.3}. 
\item The direct sum $\g\oplus \g = \g_{\text{diag}}\oplus \g_R$ with  
     the metric $g_2((X,Y),(X,Y)) = g(X,X) - g(Y,Y)$ for  
     $(X,Y)\in \g\oplus \g$, 
     where $\g_{\text{diag}}=\{(X,X):X\in\g\}$ is isomorphic to $\g$,  
     and where the subalgebra $\g_R=\{((R+1)X,(R-1)X):X\in\g\}$ is  
     isomorphic to the Lie algebra $(\g,b_R)$ with bracket  
     $b_R(X,Y)=[RX,Y]+[X,RY]$, which  
     again is isomorphic to $(\g^*,b')$, see \nmb!{2.8}. 
\endroster 
\endproclaim 
 
\demo{Proof} 
For an $R$-matrix $R$ the mapping   
$(R+1,R-1):(\g,b_R)\to \g\x \g$ is a homomorphism of Lie algebras  
into the direct product by lemma \nmb!{2.10}, which is injective.  
Also by lemma \nmb!{2.8} the   
mapping $g:(\g,b_R)\to(\g^*,b')$ is an isomorphism of Lie algebras.  
The direct sum Lie algebra $\g\oplus\g$ admits a decomposition into  
Lie subalgebras 
$$\align 
\g\oplus\g &= \{(X,X):X\in \g\} \oplus \{((R+1)Y,(R-1)Y):Y\in\g\},  
     \quad\text{ where}\\ 
(U,V) &= (X,X) + ((R+1)Y,(R-1)Y),\\ 
2X&= R(V-U)+V+U,\qquad 2Y=U-V, 
\endalign$$ 
which are isotropic: 
$$ 
g_2((R+1)Y,(R-1)Y) = g((R+1)Y,(R+1)Y)-g((R-1)Y,(R-1)Y) = 0 
$$ 
since  
$R$ is skew symmetric for $g$.  
\qed\enddemo 
 
\subhead\nmb.{2.23}. Remark \endsubhead 
The construction of an $R$-matrix on a semisimple metrical Lie  
algebra $(\g,g)$ reduces to the construction of a Manin decomposition  
$\g\oplus\g=\g_-\oplus\g_+$ of the metrical Lie algebra  
$(\g\oplus\g,g\oplus(-g))$ where $\g_-=\g_{\text{diag}}$ is the  
diagonal subalgebra.  
 
\head\totoc\nmb0{3}. Notation on Lie Groups \endhead 
 
\subhead\nmb.{3.1}. Notation for Lie groups 
\endsubhead 
Let $G$ be a Lie group with Lie algebra $\g=T_eG$,  
multiplication $\mu:G\x G\to G$, and for $g\in G$  
let $\mu_g, \mu^g:G\to G$ denote the left and right translation,  
$\mu(g,h)=g.h=\mu_g(h)=\mu^h(g)$.  
 
Let $L,R:\g\to \X(G)$ be the left  
and right invariant vector field mappings, given by  
$L_X(g)=T_e(\mu_g).X$ and $R_X=T_e(\mu^g).X$, respectively.  
They are related by $L_X(g)=R_{\Ad(g)X}(g)$. 
Their flows are given by  
$$\Fl^{L_X}_t(g)= g.\exp(tX)=\mu^{\exp(tX)}(g),\quad 
\Fl^{R_X}_t(g)= \exp(tX).g=\mu_{\exp(tX)}(g).$$ 
 
Let $\ka^l,\ka^r:\in\Om^1(G,\g)$ be the left and right Maurer-Cartan  
forms, given by $\ka^l_g(\xi)=T_g(\mu_{g\i}).\xi$ and  
$\ka^r_g(\xi)=T_g(\mu^{g\i}).\xi$, respectively. These are the  
inverses to $L,R$ in the following sense: $L_g\i=\ka^l_g:T_gG\to\g$  
and $R_g\i=\ka^r_g:T_gG\to\g$. They are related by  
$\ka^r_g=\Ad(g)\ka^l_g:T_gG\to\g$ and they satisfy the Maurer-Cartan  
equations $d\ka^l+\frac12[\ka^l,\ka^l]^\wedge =0$ and  
$d\ka^r-\frac12[\ka^r,\ka^r]^\wedge =0$.  
 
The (exterior) derivative of the function $\Ad:G\to GL(\g)$ can be  
expressed by 
$$d\Ad = \Ad.(\ad\o\ka^l) = (\ad\o \ka^r).\Ad,$$ 
which follows from 
$d\Ad(T\mu_g.X) = \frac d{dt}|_0 \Ad(g.\exp(tX)) 
= \Ad(g).\ad(\ka^l(T\mu_g.X))$. 
 
\subhead\nmb.{3.2}. Analysis on Lie groups \endsubhead 
Let $V$ be a vector space. For $f\in C^\infty(G,V)$ we have  
$df\in\Om^1(G;V)$, a 1-form on $G$ with values in $V$. We define the  
\idx{\it left derivative} $\de f=\de^lf:G\to L(\g,V)$ of $f$ by  
$$ 
\de f(x).X:= df.T_e(\mu_x).X = (L_Xf)(x)\text{ for }x\in G, X\in \g. 
$$ 
 
\proclaim{Result} \cit!{27} 
\roster 
\item For $f\in C^\infty(G,\Bbb R)$ and $g\in C^\infty(G,V)$ we have  
      $\de(f.g)=f.\de g + \de f\otimes g$, where we use  
      $\g^*\otimes V\to L(\g,V)$. 
\item For $f\in C^\infty(G, V)$ we have  
      $\de\de f(x)(X,Y)-\de\de f(x)(Y,X)=\de f(x)([X,Y])$. 
\item Fundamental theorem of calculus: 
      For $f\in C^\infty(G,V)$, $x\in G$, $X\in \g$ we have  
$$ 
f(x.\exp(X)) - f(x) = \Bigl(\int_0^1 \de f(x.\exp(tX))\,dt\Bigr)(X). 
$$ 
\item Taylor expansion with remainder: 
      For $f\in C^\infty(G,V)$, $x\in G$, $X\in \g$ we have  
$$ 
f(x.exp(X)) = \sum_{j=0}^N \frac1{j!}\de^jf(x)(X^j) +  
\int_0^1\frac{(1-t)^N}{N!}\de^{N+1}f(x.\exp(tX))\,dt\;(X^{N+1}). 
$$ 
\item For $f\in C^\infty(G,V)$ and $x\in G$ the formal Taylor series 
$$ 
\operatorname{Tay}_xf =  
     \sum_{j=0}^\infty \frac1{j!}\de^jf(x):\bigotimes \g \to \Bbb R  
$$ 
      factors to a linear functional on the universal enveloping  
      algebra: $\Cal U(\g)\to \Bbb R$. If for $A\in \Cal U(\g)$ we  
      denote by $L_A$ the associated left invariant differential  
      operator on $G$, we have  
      $\langle A, \operatorname{Tay}_xf \rangle= (L_Af)(x)$ 
\endroster 
\endproclaim 
 
\subhead\nmb.{3.3}. Vector fields and differential forms \endsubhead 
For $f\in C^\infty(G,\g)$ we get a smooth vector field $L_f\in \X(G)$  
by $L_f(x):=T_e(\mu_x).f(x)$. This describes an isomorphism  
$L:C^\infty(G,\g)\to \X(G)$. 
If $h\in C^\infty(G,V)$ then we have 
$L_fh(x)= dh(L_f(x)) = dh.T_e(\mu_x).f(x) = \de h(x).f(x)$, for which  
we write shortly $L_fh=\de h.f$. 
 
For $g\in C^\infty(G,{\tsize\bigwedge}^k\g^*)$ we get a k-form  
$L_g\in\Om^k(G)$ by the prescription  
$(L_g)_x = g(x)\o {\tsize\bigwedge}^kT_x(\mu_{x\i})$. This gives an isomorphism  
$L:C^\infty(G,{\tsize\bigwedge}\g) \to \Om(G)$.  
 
\proclaim{Result} \cit!{27} 
\roster 
\item For $f,g\in C^\infty(G,\g)$ we have  
$$ 
[L_f,L_g]_{\X(G)} = L_{K(f,g)}, 
$$ 
     where  
     $K(f,g)(x):=[f(x),g(x)]_{\g} + \de g(x).f(x)-\de f(x).g(x)$, 
     or shorter  
     $K(f,g)=[f,g]_{\g}+ \de g.f- \de f.g$. 
\item For $g\in C^\infty(G,{\tsize\bigwedge}^k\g^*)$ and $f_i\in C^\infty(G,\g)$ we  
       have $L_g(L_{f_1},\dots,L_{f_k})=g.(f_1,\dots,f_k)$. 
\item For $g\in C^\infty(G,{\tsize\bigwedge}^k\g^*)$ the exterior derivative is  
       given by  
$$ 
d(L_g) = L_{\de^\wedge g + \partial^{\g}\o g}, 
$$ 
     where $\de^\wedge g:G\to {\tsize\bigwedge}^{k+1}\g^*$ is given by  
$$\de^\wedge g(x)(X_0,\dots,X_k) =  
     \sum_{i=0}^k(-1)^i\de g(x)(X_i)(X_0,\dots,\widehat{X_i},\dots,X_k),$$ 
     and where $\partial^{\g}$ is the Chevalley differential on  
     ${\tsize\bigwedge}\g^*$. 
\item For $g\in C^\infty(G,{\tsize\bigwedge}^k\g^*)$ and $f\in C^\infty(G,\g)$ the  
       Lie derivative is given by  
$$ 
\L_{L_f}L_g = L_{\L^{\g}_f\o g + \L^\de_f g}, 
$$ 
where  
$$\align 
(\L^{\g}_f g)(x)(X_1,\dots,X_k) &= \sum_{i}(-1)^i  
g(x)([f(x),X_i],X_1,\dots,\widehat{X_i},\dots, X_k),\\ 
(\L^{\de}_f g)(x)(X_1,\dots,X_k)  
     &= \de g(x)(f(x))(X_1,\dots,X_k) +\\ 
&+\sum_{i}(-1)^i  
     g(x)(\de f(x)(X_i),X_1,\dots,\widehat{X_i},\dots, X_k). 
\endalign$$ 
\endroster 
\endproclaim 
 
\subhead\nmb.{3.4}. Multi vector fields and the Schouten-Nijenhuis  
bracket \endsubhead  
Recall that on a manifold $M$ the space of multi vector fields  
$\Ga({\tsize\bigwedge} TM)$ carries the Schouten-Nijenhuis bracket, given by  
$$\align 
[X_1\wedge &\dots\wedge X_p,Y_1\wedge \dots\wedge Y_q] = \tag1\\ 
&= \sum_{i,j}(-1)^{i+j}[X_i,Y_j]  
     \wedge \dots\widehat{X_i}\dots \wedge X_p\wedge  
     Y_1\wedge \dots\widehat{Y_j}\dots\wedge Y_q.   
\endalign$$ 
See \cit!{28} for a presentation along the lines used here.  
This bracket has the following properties: Let $U\in\Ga({\tsize\bigwedge}^uTM)$,  
$V\in\Ga({\tsize\bigwedge}^vTM)$, $W\in\Ga({\tsize\bigwedge}^wTM)$, and $f\in C^\infty(M,\Bbb R)$.  
Then 
$$\align 
[U,V] &= -(-1)^{(u-1)(v-1)}[V,U] \\ 
[U,[V,W]] &= [[U,V],W] 
     + (-1)^{(u-1)(v-1)}[V,[U,W]]\\ 
[U,V\wedge W] &= [U,V]\wedge W  
     + (-1)^{(u-1)v} V\wedge [U,W]\\ 
[f,U] &= -\bar \imath (df)U, 
\endalign$$ 
where $\bar \imath(df)$ is the insertion operator  
${\tsize\bigwedge}^kTM\to {\tsize\bigwedge}^{k-1}TM$, the adjoint of  
$df\wedge(\quad):{\tsize\bigwedge}^lT^*M\to{\tsize\bigwedge}^{l+1}T^*M$.  
 
For a Lie group $G$ we have an isomorphism 
$L:C^\infty(G,{\tsize\bigwedge} \g)\to \Ga({\tsize\bigwedge} TG)$ 
which is given by   
$L(u)_x = {\tsize\bigwedge} T(\mu_x).u(x)$, via left trivialization.  
For $u\in C^\infty(G,{\tsize\bigwedge}^u\g)$ we have  
$\de u:G\to L(\g,{\tsize\bigwedge}^u\g)=\g^*\otimes {\tsize\bigwedge}^u\g$, 
and with respect to   
the one component in $\g^*$ we can consider the insertion operator 
$\bar\imath(\de u(x)):{\tsize\bigwedge}^k\g\to {\tsize\bigwedge}^{k+u}\g$. 
In more detail: if   
$u=f.U$ for $f\in C^\infty(G,\Bbb R)$ and 
$U\in {\tsize\bigwedge}^u\g$, then we put  
$\bar\imath(\de f(x).U)V= U\wedge \bar\imath(\de f(x))(V)$. 
 
For the Lie algebra $\g$ we also have the algebraic  
Schouten-Nijenhuis bracket  
$[\quad,\quad]^\g:{\tsize\bigwedge}^p\g \x {\tsize\bigwedge}^q\g \to 
     {\tsize\bigwedge}^{p+q-1}\g$  
which is given by formula \thetag1, applied to this purely algebraic  
situation.   
 
\proclaim{Proposition} 
For $u\in C^\infty(G,{\tsize\bigwedge}^u\g)$ and  
$v\in C^\infty(G,{\tsize\bigwedge}^v\g)$ the   
Schouten-Nijenhuis bracket is given by  
$$ 
[L(u),L(v)] = L([u,v]^\g - \bar\imath(\de u)(v) 
     + (-1)^{(u-1)(v-1)}\bar\imath(\de v)(u)).\tag2   
$$ 
\endproclaim 
 
\demo{Proof} 
This follows from formula \thetag1 applied to  
$$ 
[L(f.X_1\wedge \dots\wedge X_p), 
     L(g.Y_1\wedge \dots\wedge Y_q)], 
$$ 
where $f,g\in C^\infty(G,\Bbb R)$ and $X_i,Y_j\in \g$, and then by  
applying \nmb!{3.3}.\therosteritem1.  
\qed\enddemo 
 
\head\totoc\nmb0{4}.  
Lie Poisson Groups and Double Groups  
\endhead 
 
\subhead\nmb.{4.1}. Lie Poisson groups \endsubhead 
A \idx{\it Poisson structure} on a Lie group is a tensor field  
$\La\in \Ga({\tsize\bigwedge}^2TG)$ such that $\{f,g\}:=\langle df\wedge dg,\La  
\rangle$ defines a Lie bracket on $C^\infty(G,\Bbb R)$.  
If we let $\La=L(\la)$ for $\la\in C^\infty(G,{\tsize\bigwedge}^2\g)$ in the notation of  
\nmb!{3.4}, then $\La$ is a Poisson structure if and only if for the  
Schouten bracket we have  
$[\La,\La]=0$. By proposition \nmb!{3.4} this is  
equivalent to   
$$ 
[\la(g),\la(g)]^{\g} = 2\bar\imath(\de \la(g))(\la(g))\quad 
\text{ for all }g\in G.\tag1 
$$ 
 
A \idx{\it Lie-Poisson} group \cit!{11} is a Lie group $G$  
together with a Poisson structure $\La\in \Ga({\tsize\bigwedge}^2TG)$ such that the  
multiplication $\mu:G\x G\to G$ is a Poisson map, i.e\. the pull back  
mapping $\mu^*:C^\infty(G,\Bbb R)\to C^\infty(G\x G,\Bbb R)$ is a  
homomorphism for the Poisson brackets.  
This is equivalent to 
any of the following properties \therosteritem2 -- \therosteritem6  
for $p=2$ (see \cit!{21}).  
Such a 2-vector field $\La$ is also called a  
\idx{\it Lie-Poisson structure}. 
 
\proclaim{Lemma} 
For $\La\in\Ga({\tsize\bigwedge}^pTG)$ the following assertions  
\therosteritem2--\therosteritem6 are equivalent: 
\roster 
\item [2] $\La$ is multiplicative in the sense that 
$$ 
\La(gh) = {\tsize\bigwedge}^pT(\mu_g).\La(h) +  
     {\tsize\bigwedge}^pT(\mu^h).\La(g)   
     \text{ for all }g,h\in G.  
$$ 
\item (assuming that $G$ is connected) $\La(e)=0$ and the Schouten  
       bracket $\L_{L_X}\La=[L_X,\La]$ is  
       left invariant for each left invariant vector field $L_X$ on  
       $G$. 
\item (assuming that $G$ is connected) $\La(e)=0$ and the Schouten  
       bracket $\L_{R_X}\La=[R_X,\La]$ is  
       right invariant for each right invariant vector field $R_X$ on  
       $G$. 
\item If we let $\La=L(\la)$ for  
       $\la\in C^\infty(G,{\tsize\bigwedge}^p\ga)$ in the   
       notation of \nmb!{3.4}, then  
$$ 
\la(gh) = {\tsize\bigwedge}^p\Ad(h\i).\la(g) +  
     \la(h)\quad\text{ for all }g,h\in G.  
$$ 
       This has the following meaning: Consider the right semidirect  
       product $G\ltimes {\tsize\bigwedge}^p \g$ with multiplication  
       $(x,U).(y,V)=(xy,\Ad(y\i)U + V)$. Then the above equation  
       holds if and only if $x\mapsto (x,\la(x))$ is a homomorphism of  
       Lie groups.  
\item  $\La:G\to {\tsize\bigwedge}^pTG$ is a homomorphism of Lie groups,  
       where $L:G\ltimes {\tsize\bigwedge}^p\g\cong {\tsize\bigwedge}^pTG$.   
\endroster 
A Poisson structure $\La$ on $G$ is a Lie-Poisson structure if and  
only if these conditions \therosteritem2--\therosteritem6 are  
satisfied for $p=2$. 
\endproclaim 
 
\demo{Proof}  
For the proof of the equivalence of conditions  
\therosteritem2--\therosteritem4 see \cit!{21}, the equivalence to  
\therosteritem5 and \therosteritem6 is obvious. 
 
We prove the last assertion. It follows from 
$$\align 
\{\mu^*f,\mu^*g\}_{G\x G}(x,y)  
&= \langle d(f\o\mu)\wedge d(g\o\mu),\La(x)\otimes \La(y) \rangle\\ 
&= (df(xy)\wedge dg(xy)).{\tsize\bigwedge}^2T_{(x,y)}\mu.(\La(x),\La(y))\\ 
&= (df(xy)\wedge dg(xy)). 
     {\tsize\bigwedge}^2(T_y(\mu_x)+T_x(\mu^y)).(\La(x),\La(y))\\  
&= (df(xy)\wedge dg(xy)).({\tsize\bigwedge}^2T_y(\mu_x)\La(y) 
     +{\tsize\bigwedge}^2T_x(\mu^y)\La(x))\\ 
\endalign$$ 
compared with 
$$ 
(\mu^*\{f,g\}_G)(x,y) = (df\wedge dg).\La(xy). \qed 
$$ 
\enddemo 
 
Note that if $\La_1:G\to {\tsize\bigwedge}^{p_1}TG$ and  
$\La_2:G\to {\tsize\bigwedge}^{p_2}TG$ are homomorphisms of groups with  
$\pi\o\La_i = \Id_G$, then their Schouten bracket  
$[\La_1,\La_2]:G\to {\tsize\bigwedge}^{p_1+p_2-1}TG$ 
has the same property. This follows from \cit!{21} and the  
equivalence to \therosteritem6 from above.

\proclaim{\nmb.{4.2}. Theorem} \cit!{11} 
If $(G,\La)$ is a Lie-Poisson group with Lie algebra $\g$ then by  
$b':\g\to{\tsize\bigwedge}^2\g$ we get a Lie bialgebra structure on  
$\g$, where 
$b'(X)=(\L_{L_X}\La)(e)=\de \la(e)X$, where $\L$ denotes the Lie  
derivative. 
 
If $(\g,b,b')$ is a Lie bialgebra and $G$ is a simply connected  
Lie group associated to $\g$, then the cocycle $b'$ integrates to a  
unique Lie Poisson structure $\La\in\Ga({\tsize\bigwedge}^2TG)$ on $G$. 
\endproclaim 
 
\demo{Proof} See \cit!{11} and \cit!{21} for other proofs.  
By conditions \nmb!{4.1}.\therosteritem5 and \therosteritem6 any  
multiplicative 2-vector-field $\La$ 
is a homomorphism of Lie-groups  
$$\CD 
G @>{\La}>> {\tsize\bigwedge}^2TG\\ 
@|          @A{\cong}AA\\ 
G @>{(\Id,\la)}>> G \ltimes{\tsize\bigwedge}^2\g 
\endCD$$ 
and the induced Lie algebra homomorphism then is  
$$\align 
T_e(\La).X:=&(X,\L_{L_X}\La(e)) \\ 
=& (X,\de\la(e).X) \quad\text{ by Proposition \nmb!{3.4}.\thetag2}\\ 
=& (\Id_{\g},b')(X), 
\endalign$$ 
and conversely any 2-cocycle $b':\g\to{\tsize\bigwedge}^2\g$ integrates  
to a Lie group homomorphism if $G$ is supposed to be simply  
connected.  
 
It remains to show that $b':{\tsize\bigwedge}^2\g^* \to \g^*$ satisfies the  
Jacobi identity if and only if \nmb!{4.1}.\therosteritem1 holds.  
Let us take the left derivative $\de$ at $e$ of equation  
\nmb!{4.1}.\therosteritem1 and get 
$$\align 
0&=2[\de\la(e),\la(e)]^{\g} -  
2\bar\imath(\de^2\la(e))\la(e) - 2\bar\imath(\de\la(e))\de\la(e) \\ 
&= 0 - 0 - [\de\la(e),\de\la(e)]^{\text{NR}}, 
\endalign$$  
so that the Nijenhuis-Richardson bracket of  
$b'=\de\la(e):{\tsize\bigwedge}^2\g^*\to\g^*$ with itself vanishes.  
This just means that $b'$ is a Lie bracket on $\g^*$, see 
\cit!{30}. 
 
For the converse note first that if $\La:G\to{\tsize\bigwedge}^2TG$ is a  
homomorphism of Lie groups then also the Schouten bracket  
$[\La,\La]:G\to{\tsize\bigwedge}^3TG$   
is a homomorphism of Lie groups.  
But if  
$b'=\de\la(e)$ is a Lie bracket on $\g^*$ then the computation above  
shows that  
$\de\Bigl([\la,\la]^{\g}-2\bar\imath(\de\la)\la\Bigr)(e)=0$ 
so that the associated Lie algebra homomorphism is just  
$(\Id,0):\g\to\g\ltimes {\tsize\bigwedge}^3\g$. But then  
$[\La,\La]=0$.  
\qed\enddemo

\subhead\nmb.{4.3}. Affine Poisson structures 
\endsubhead 
An affine Poisson structure on a Lie group $G$ is a Poisson structure  
$\La$ such that $\La_l$ is a Lie Poisson structure or equivalently $\La_r$  
is a Lie Poisson structure, where  
$$\alignat2 
\La_l(g) &= \La(g) - T(\mu_g)\La(e),  
     &\quad \La_l &= \La - L_{\La(e)},\tag1\\ 
\La_r(g) &= \La(g) - T(\mu^g)\La(e),  
     &\quad \La_r &= \La - R_{\La(e)}.\tag2\\ 
\endalignat$$ 
For a Poisson structure $\La$ we also have  
$$\alignat2 
\La_l &= L(\la_l), &\quad \la_l(g) &= \la(g) - \la(e),\tag{1'}\\ 
\La_r &= L(\la_r), &\quad \la_r(g) &= \la(g) - \Ad(g\i)\la(e),\tag{2'}\\ 
\endalignat$$ 
and $\La$ is an affine Poisson structure if and only if  
$$ 
\la(gh) = \Ad(h\i)\la(g) + \la(h) - \Ad(h\i)\la(e). \tag{3} 
$$

\subhead\nmb.{4.4}. Lie groups with exact Lie bialgebras\endsubhead 
Let $G$ be a Lie group with Lie algebra $\g$. 
Suppose we have a solution $C\in{\tsize\bigwedge}^2\g$ of the \thetag{mYBE},  
so that $b'=\partial C$ is a Lie bialgebra structure for $(\g,\g^*)$. 
Then we can write down explicitly the Lie-Poisson structure on any  
(even not connected) Lie group with Lie algebra $\g$, as follows.  
 
We consider $\La_{\pm}:G\to {\tsize\bigwedge}^2TG$ qiven by  
$\La_\pm(g):= T(\mu_g)C\pm T(\mu^g)C$.  
Then obviously $\La_-$ is multiplicative and $\La_+$ is affine with  
$(\La_+)_l = \La_-$ and $(\La_+)_r = -\La_-$.  
In the notation of  
\nmb!{4.1} we have $\la_\pm (g) = C\pm\Ad(g\i)C$, and  
$$ 
b'_\pm(X)=\de \la_\pm(e)X =  
\pm(\de({\tsize\bigwedge}^2(\Ad\o \operatorname{Inv}))(e)X)C = \mp\ad(X)C =  
\mp(\partial_bC)(X), 
$$  
and since $C$ satisfies \thetag{mYBE}, the tensor fields $\La_\pm$  
are Poisson structures. 
 
\subhead\nmb.{4.5}. Manin decompositions and Lie-Poisson structures 
\endsubhead  
Let $\g=\g_+\oplus\g_-$ be a Manin decomposition  
of a metrical Lie algebra $\g$, and let  
$\pr_{\pm}:\g\to\g_{\pm}$ be the corresponding projections. 
Then by \nmb!{2.20} the operator 
$R=\pr_+-\pr_-$ is a solution of \thetag{1-mYBE} $B_R+b=0$. 
 
So by \nmb!{4.4} a Manin decomposition defines a canonically associated  
Lie-Poisson structure on each (even not connected) Lie group $G$ with  
Lie algebra $\g$, as follows: 
Let $C=R\o g\i\in{\tsize\bigwedge}\g$ be the associated exact bialgebra  
structure, and    
consider $\La_{\pm}:G\to {\tsize\bigwedge}^2TG$ qiven by  
$$\La_\pm(g):= T(\mu_g)C\pm T(\mu^g)C.\tag1$$ 
Then in the notation of  
\nmb!{4.1} we have $\la_\pm (g) = C\pm\Ad(g\i)C$, and  
$b'_\pm(X)=\de \la_\pm(e)X =  
\pm{\tsize\bigwedge}^2(\de(\Ad\o \operatorname{Inv})(e)X)C = \mp\ad(X)C =  
\mp (\partial_bC)(X)$.  
The tensor field $\La_-$ is a real analytic Lie-Poisson structure and  
$\La_+$ is a real analytic affine Poisson structure with  
$(\La_+)_l = \La_-$ and $(\La_+)_r = -\La_-$.   
Since $\La_+(e)=C$ is non-degenerate, the affine Poisson structure  
$\La_+$ is non-degenerate on an open subset of $G$. If $G$ is connected  
this open subset is also dense since the real analytic Poisson  
structure cannot be degenerate on an open subset.  
 
We shall investigate this kind of structure in much more details  
below.  
 
\subhead\nmb.{4.6}. Gau\ss{} decompositions and Lie-Poisson  
structures \endsubhead 
Let $G$ be a Lie group with a metrical Lie algebra $(\g,g)$.  
 From \nmb!{2.21} we know that any solution $R$ of the $R$-matrix  
equation can be described as follows.  
There is a Gau\ss{} decomposition  
$\g=\g_+\oplus\g_0\oplus\g_-$ with $\g_\pm$ isotropic and dual to  
each other, and with $g$ non-degenerate on $\g_0$.  
Let $\pr_{\pm,0}:\g\to\g_{\pm,0}$ be the  
orthogonal projections.  
Then $R$ is of the following form: 
$$R=\pr_+\oplus (R_0\o \pr_0)\oplus(-\pr_-),\tag1$$ 
where $R_0$ is a solution of \thetag{1-mYBE} on $\g_0$ without  
eigenvalues $1$ or $-1$ (without fixed points). 
 
Let $X_i$ be a basis of $\g_+$ with $Y_i$ the dual basis of $\g_-$,  
and let $Z_j$ be an orthonormal basis of $\g_0$, all   
with respect to the inner product $g$ on $\g=\g_+\oplus\g_0\oplus\g_-$.  
Let $R_0(Z_j)=\sum_k R_j^kZ_k=\sum_k C^{kj}Z_k$ be the (skew  
symmetric) matrix representation of $R_0$ with respect to   
the basis $Z_j$.  
Then  
$$ 
\pr_+(U)=\sum_i X_i.g(U,Y_i),\quad  
\pr_0(U) = \sum_j Z_j.g(U,Z_j),\quad 
\pr_-(U)=\sum_i Y_i.g(U,X_i), 
$$ 
so that  
$$\align 
R&= \pr_+  - \pr_- + (R_0\o\pr_0)  
= \Bigl(\sum_i X_i\wedge Y_i +  
     \sum_{j,k}R^k_jZ_k\otimes Z_j\Bigr)\o g\tag2\\ 
C:&=R\o g\i = \sum X_i\wedge Y_i+\sum_{j<k}C^{jk}Z_j\wedge Z_k.  
\endalign$$ 
Let us consider  
$\La_{\pm}:G\to {\tsize\bigwedge}^2TG$ qiven by  
$$
\La_\pm(g):= T(\mu_g)C\pm T(\mu^g)C.
\tag3$$ 
Then in the notation of  
\nmb!{4.1} we have $\la_\pm (g) = C\pm\Ad(g\i)C$, and  
$$ 
b'_\pm(X)=\de \la_\pm(e)X = \pm(\de(\Ad\o \operatorname{Inv})(e)X)C  
= \mp\ad(X)C = \mp (\partial_bC)(X).$$  
Since $R$ was a solution of \thetag{1-YBE} the tensor field  
$\La_-$ is a real analytic Lie-Poisson structure and $\La_+$ is a real  
analytic affine Poisson structure with  
$(\La_+)_l = \La_-$ and $(\La_+)_r = -\La_-$.   
Since $\La_+(e)=C$ is non-degenerate, the affine Poisson structure  
$\La_+$ is non-degenerate on an open subset of $G$. If $G$ is connected  
this open subset is also dense since the real analytic Poisson  
structure cannot be degenerate on an open subset.  
  
\head\totoc\nmb0{5}. Explicit Formulas for Poisson Structures on  
Double Lie Groups \endhead  
 
\subhead\nmb.{5.1}. The setting \endsubhead 
It turns out that in the situation of \nmb!{4.5} one can get very  
useful explicit formulae. Let us explain this setting once more,  
which will be used for the rest of this paper.   
 
Let $G$ be any Lie group with a metrical Lie algebra $(\g,\ga)$  
and suppose that it admits a Manin decomposition  
$(\g=\g_+\oplus\g_-,\ga)$. 
Let $\pr_{\pm}:\g\to\g_{\pm}$ be the corresponding  
projections. By \nmb!{2.20} the operator  
$R=\pr_+-\pr_-$ is a solution of \thetag{1-mYBE} $B_R+b=0$. 
 
\subhead Simplified notation \endsubhead 
In order to compactify the notation we will use the following  
shorthand, in the rest of this paper: For  
$U\in\bigotimes^p\g$ etc\. and for $g\in G$ we let 
$$ 
gU = g.U = \bigotimes^pT(\mu_g)U,  
     \qquad Ug = U.g = \bigotimes^pT(\mu^g)U. 
$$ 
 
\bigskip 
Let $X_i$ be a basis of $\g_+$ with $Y_i$ the dual basis of $\g_-$  
with respect to the inner product $\ga$ on $\g=\g_+\oplus\g_-$. Then  
$$\align 
\pr_+(Z)&=\sum_i \ga(Z,Y_i).X_i=(\sum_i Y_i\otimes X_i)\ga(Z),\\ 
\pr_-(Z)&=\sum_i \ga(Z,X_i)Y_i=(\sum_i X_i\otimes Y_i)\ga(Z),\quad  
     \text{ so that}\\ 
\pr_+&=\Bigl(\sum_iY_i\otimes X_i\Bigr)\o \ga = C_+\o\ga,\\ 
\pr_-&=\Bigl(\sum_iX_i\otimes Y_i\Bigr)\o \ga = C_-\o\ga,\quad  
     \text{ where}\\ 
C_+&=\sum_iY_i\otimes X_i, \quad C_-=\sum_iX_i\otimes Y_i. 
\endalign$$ 
Then we have  
$$\align 
R&= \pr_+ - \pr_- = (\sum_i Y_i\wedge X_i)\o \ga,\tag1\\ 
C&= R\o\ga\i = C_+-C_-= \sum Y_i\wedge X_i.  
\endalign$$ 
Then we  
consider $\La_{\pm}:G\to {\tsize\bigwedge}^2TG$ qiven by (note the factor  
$\frac12$) 
$$
\La_\pm(g):= \frac12(gC\pm Cg).
\tag2$$ 
Then in the notation of  
\nmb!{4.1} we have $\la_\pm (g) = \frac12(C\pm\Ad(g\i)C)$, and  
$b'_\pm(X)=\de \la_\pm(e)X =  
\pm\tfrac12{\tsize\bigwedge}^2(\de(\Ad\o \operatorname{Inv})(e)X)C  
= \mp\tfrac12\ad(X)C = \mp\tfrac12(\partial_bC)(X)$.  
The tensor field $\La_-$ is a real analytic Lie-Poisson structure and  
$\La_+$ is a real analytic affine Poisson structure with  
$(\La_+)_l = \La_-$ and $(\La_+)_r = -\La_-$.   
Since $\La_+(e)=C$ is non-degenerate, the affine Poisson structure  
$\La_+$ is non-degenerate on an open subset of $G$. If $G$ is connected  
this open subset is also dense since the real analytic Poisson  
structure cannot vanish on an open subset.

\proclaim{\nmb.{5.2}. Lemma} In the setting of \nmb!{5.1} we have: 
$$\align 
\La_+(a) &= aC_+ - C_-a= C_+a - aC_- \tag1\\ 
&= \sum_i \Bigl(aY_i\otimes aX_i - X_ia\otimes Y_ia\Bigr) 
     = \sum_i \Bigl(Y_ia\otimes X_ia - aX_i\otimes aY_i\Bigr)\\ 
\La_-(a)  
&= aC_+ - C_+a = C_-a - aC_-\tag2\\ 
&= \sum_i \Bigl(aY_i\otimes aX_i-Y_ia\otimes X_ia\Bigr) 
     = \sum_i \Bigl(X_ia\otimes Y_ia - aX_i\otimes aY_i\Bigr)\\  
\endalign$$ 
\endproclaim 
 
\demo{Proof} 
The tensor fields do not look skew symmetric but observe that  
$$ 
aC_+ + aC_- = C_+a + C_-a. \tag3 
$$ 
This is equivalent to  
$C_+ + C_- = \bigotimes^2\Ad(a\i)(C_+ + C_-)$ 
which, when composed with $\ga$, in $L(\g,\g)$ just says that 
$\Id_\g = \pr_+ + \pr_- = \Ad(a\i)\Id_\g\Ad(a)$. 
Using \thetag5 we have  
$$\align 
\La_+(a) &= \tfrac12(aC+Ca) = \tfrac12(aC_+ - aC_- + C_+a - C_-a)\\ 
&= C_+a - aC_- = aC_+ - C_-a,\\ 
\La_-(a) &= \tfrac12(aC-Ca) = \tfrac12(aC_+ - aC_- - C_+a + C_-a)\\ 
&= aC_+ - C_+a = C_-a - aC_-.\qed 
\endalign$$ 
\enddemo 
 
\subhead\nmb.{5.3}. The subgroups and the Poisson structures  
\endsubhead  
In the setting of \nmb!{5.1} we consider now the Lie subgroups  
$G_\pm\subset G$ corresponding to the isotropic Lie subalgebras 
$\g_\pm$, and we consider the mappings 
$$\alignat2 
\ph: &G_+\x G_- \to G, &\quad  \ph(g,u):&= g.u \in G,\\ 
\ps: &G_-\x G_+ \to G, &\quad  \ps(v,h):&= v.h \in G. 
\endalignat$$  
Both are diffeomorphisms on open neighbourhoods of $(e,e)$. 
We will use $g,u$ and $v,h$ as local `coordinates' near $e$. 
So, we have, at least locally in an open neighborhood $U$ of $e$ in  
$G$, well defined projections $p^+_l,p^+_r:G\supset U\to G_+$ and  
$p^-_l,p^-_r:G\supset U\to G_-$ which play the role of momentum mappings: 
$$\alignat2 
p^+_l(g.u) :&= g, &\quad  p^+_r(v.h) :&= h \in G_+,\\ 
p^-_l(v.h) :&= v, &\quad  p^-_r(g.u) :&= u \in G_-. 
\endalignat$$  
When $\ph$ (or equivalently $\ps$) is a global diffeomorphism (this  
is consistent for simply connected $G$ with completeness of the the 
dressing vector fields; in these cases we will call $G$ a 
\idx{\it complete double group})   
then the mappings $p^\pm_{l,r}$ are globally defined.  
 
\proclaim{Remark} If the subgroup $G_+$ is compact then the double  
group $G$ is complete. Similarly for $G_-$. 
\endproclaim 
Indeed, there exists a $G$-invariant Riemann metric on the  
homogeneous space $G/G_+$. Then $G$ acts on $G/G_+$ by isometries  
locally transitively, hence transitively. This means that $G=G_+.G_-$  
globally and that $G_+\cap G_-$ is finite. 
 
\proclaim{\nmb.{5.4}. Theorem} 
In the setting above,   
the following tensor fields are Lie-Poisson structures on  
the group $G_+$ and $G_-$, respectively, corresponding to the Lie  
bialgebra structures on $\g_+$ and $\g_-$ induced from the Manin  
decomposition: 
$$\align 
\La^{G_+}(g) &= g((\Id_\g\otimes\Ad^G(g\i)\pr_+\Ad^G(g))C_-) 
     \in {\tsize\bigwedge}^2TG_+\tag1 \\ 
&= g(-(\Ad^G(g\i)\otimes\pr_+\Ad^G(g\i))C_-)\\  
&= {\tsize\sum}_i gX_i \otimes \pr_+(\Ad^G(g)Y_i)g\\ 
&= -{\tsize\sum}_i X_ig \wedge g\pr_+(\Ad^G(g\i)Y_i),\\ 
\La^{G_-}(u) &= u((\Id_\g\wedge \Ad^G(u\i)\pr_-\Ad^G(u))C_+) 
     \in{\tsize\bigwedge}^2TG_-\tag2\\ 
&= u(-(\Ad^G(u\i)\otimes \pr_-\Ad^G(u\i))C_+)\\ 
&= {\tsize\sum}_i uY_i \otimes \pr_-(\Ad^G(u)X_i)u\\ 
&= -{\tsize\sum}_i Y_iu \otimes u\pr_-(\Ad^G(u\i)X_i). 
\endalign$$ 
The following tensor fields are non-degenerate Poisson structures 
on the groups $G_+\x G_-$ and $G_-\x G_+$, respectively.
$$\align 
\La^\ph_+(g,u) &= \La^{G_+}(g) + \La^{G_-}(u)  
     + {\tsize\sum}_i Y_iu\wedge gX_i  
     \in{\tsize\bigwedge}^2T(G_+\x G_-),\tag3\\ 
\La^\ps_+(v,h) &= -\La^{G_+}(h) - \La^{G_-}(v)  
     + {\tsize\sum}_i vY_i\wedge X_ih  
     \in{\tsize\bigwedge}^2T(G_-\x G^+).\tag4\\ 
\endalign$$ 
Moreover they are related to the affine Poisson structures on $G$,  
i.e\. we have  
$$ 
{\tsize\bigwedge}^2 T\ph.\La^\ph_+ = \La_+ \o \ph, \qquad 
{\tsize\bigwedge}^2 T\ps.\La^\ps_+ = \La_+ \o \ps.  \tag5 
$$ 
The following tensor fields are Lie Poisson structures on the groups  
$G_+\x G_-$ and $G_-\x G_+$, respectively: 
$$\align 
\La^\ph_-(g,u) &= -\La^{G_+}(g) + \La^{G_-}(u)  
     \in{\tsize\bigwedge}^2T(G_+\x G_-),\tag6\\ 
\La^\ps_-(v,h) &= -\La^{G_+}(h) + \La^{G_-}(v) 
     \in{\tsize\bigwedge}^2T(G_-\x G_+).\tag7\\ 
\endalign$$ 
Moreover they are related to the Lie Poisson structure on $G$ which  
corresponds to $C$, i.e\.  
we have 
$$ 
T\ph. \La^\ph_- =\La_-\o \ph,\qquad T\ps. \La^\ps_-  
     =\La_-\o \ps.\tag8 
$$ 
\endproclaim 
 
\demo{Proof} 
Using \nmb!{5.2}.\thetag1 we have  
$$\align 
\La_+(gu) &= {\tsize\sum}_i (guY_i\otimes guX_i- X_igu\otimes Y_igu) \\ 
&= {\tsize\sum}_i g\Bigl(\Ad(u)Y_i\otimes \Ad(u)X_i 
     - \Ad(g\i)X_i\otimes \Ad(g\i)Y_i\Bigr)u \\ 
&= {\tsize\sum}_i g\Bigl(\Ad(u)Y_i\otimes\pr_-(\Ad(u)X_i)  
     +\Ad(u)Y_i\otimes\pr_+(\Ad(u)X_i)-\\  
&\qquad- \Ad(g\i)X_i\otimes \pr_-(\Ad(g\i)Y_i) 
     -\Ad(g\i)X_i\otimes\pr_+(\Ad(g\i)Y_i)\Bigr)u.  
\endalign$$ 
In $L(\g,\g)$ we have (compare with \nmb!{5.1}.\thetag1) 
$$\align 
\Bigl({\tsize\sum}_i \Ad(u)Y_i\otimes&\pr_+(\Ad(u)X_i)\Bigr)\o \ga 
= \pr_+\o\Ad(u)\o\pr_+\o\Ad(u\i)\\ 
&= \pr_+\o\Ad(u)\o(\Id_\g-\pr_-)\o\Ad(u\i) \\  
&= \pr_+ - \pr_+\o\Ad(u)\o\pr_-\o\Ad(u\i) = \pr_+ - 0,  
\endalign$$ 
for $\pr_+\o\Ad(u)\o\pr_- =0$ since $u\in G_-$.  
Thus we get  
$$ 
{\tsize\sum}_i \Ad(u)Y_i\otimes\pr_+(\Ad(u)X_i)  
     = {\tsize\sum}_iY_i\otimes X_i 
$$ 
and similarly we obtain   
$${\tsize\sum}_i \Ad(g\i)X_i\otimes\pr_-(\Ad(g\i)Y_i)  
     = {\tsize\sum}_iX_i\otimes Y_i,$$ 
so that  
$$\align 
\La_+(gu)  
&= g\Bigl({\tsize\sum}_i uY_i\otimes\pr_-(\Ad(u)X_i)u\Bigr)  
     +g\Bigl({\tsize\sum}_i Y_i\wedge X_i\Bigr) -\\ 
&\qquad-\Bigl({\tsize\sum}_i X_ig\otimes g\pr_+(\Ad(g\i)Y_i)\Bigr)u\\ 
&= T_{(g,u)}\ph\Bigl(\La_{G_-}(u) + \La_{G_+}(g)  
     +{\tsize\sum}_i Y_iu\wedge gX_i \Bigr), 
\endalign$$ 
which proves \thetag3 and part of \thetag5.  
In a similar way one proves \thetag4 and the other part of \thetag5. 
 
Next we check that the two expressions for $\La^{G_+}$ in \thetag1  
are the same. We have to show that the following expression vanishes 
$$\align 
&{\tsize\sum}_igX_i \otimes \pr_+(\Ad(g)Y_i)g 
+ {\tsize\sum}_i X_ig \otimes g\pr_+(\Ad(g\i)Y_i) =\\ 
&= g\Bigl({\tsize\sum}_i X_i \otimes \Ad(g\i)\pr_+(\Ad(g)Y_i) 
+ {\tsize\sum}_i \Ad(g\i)X_i \otimes \pr_+(\Ad(g\i)Y_i))\Bigr).\\ 
\endalign$$ 
The term in brackets, composed with $\ga$ from the right, is the  
following endomorphism of $\g$: 
$$\align 
\Ad&(g\i)\pr_+\Ad(g)\pr_- + \pr_+\Ad(g\i)\pr_-\Ad(g) \\ 
&=\Ad(g\i)\pr_+\Ad(g)(\Id-\pr_+) + \pr_+\Ad(g\i)(\Id-\pr_+)\Ad(g) \\ 
&=\Ad(g\i)\pr_+\Ad(g) - \Ad(g\i)\pr_+\Ad(g)\pr_+ + \pr_+  
 - \pr_+\Ad(g\i)\pr_+\Ad(g) \\ 
&=\Ad(g\i)\pr_+\Ad(g) - \pr_+ + \pr_+  
 - \Ad(g\i)\pr_+\Ad(g) =0, 
\endalign$$ 
since $\Ad(g\i)\g_+\subset \g_+$ and $\pr_+|\g_+=\Id$. 
In the same way one shows that the the two expressions for $\La^{G_-}$ in  
\thetag2 coincide, and similar computations show that all expressions  
in \thetag1 and \thetag2 are indeed skew-symmetric (which is clear  
from the beginning). 
 
Next we show that $\La^{G_+}$ is multiplicative. We have the  
following chain of equivalent assertions:  
$$\align 
\La^{G_+}(gh) &= g\La^{G_+}(h) + \La^{G_+}(g)h,\\ 
(gh)\i\La^{G_+}(gh) &= h\i\La^{G_+}(h) + h\i g\i\La^{G_+}(g)h,\\ 
{\tsize\sum}_i X_i \otimes \Ad(gh)\i\pr_+(\Ad(gh)Y_i) 
     &= {\tsize\sum}_i X_i \otimes \Ad(h\i)\pr_+(\Ad(h)Y_i) +\\ 
&\quad+{\tsize\sum}_i \Ad(h\i)X_i \otimes \Ad(gh)\i\pr_+(\Ad(g)Y_i),\\ 
\Ad(gh)\i\pr_+\Ad(gh)\pr_- &= \Ad(h\i)\pr_+\Ad(h)\pr_-\\ 
&\quad+\Ad(gh)\i\pr_+\Ad(g)\pr_-\Ad(h). \\ 
\endalign$$  
Both sides of the last equation vanish when applied to elements of  
$\g_+$, and on elements of $\g_-$ we may delete the rightmost  
$\pr_-$, so this is equivalent to  
$$\align 
\pr_+\Ad(gh) &= \Ad(g)\pr_+\Ad(h) + \pr_+\Ad(g)\pr_-\Ad(h) \\ 
&= \Ad(g)\pr_+\Ad(h) + \pr_+\Ad(g)(\Id-\pr_+)\Ad(h) \\ 
&= \Ad(g)\pr_+\Ad(h) + \pr_+\Ad(gh) - \pr_+\Ad(g)\pr_+\Ad(h), \\ 
\endalign$$ 
which is true since $\Ad(g)(\g_+)\subset \g_+$.  
 
Finally we show that the group homomorphism  
$\La^{G_+}:G_+ \to {\tsize\bigwedge}^2TG_+$ is  
associated to the bialgebra structure given by the Lie bracket on 
$\g_- @>\g>> (\g_+)^*$. 
For that we consider, as explained in \nmb!{4.1} and in the proof of  
\nmb!{4.2}:  
$$\align 
\la^{G_+}(g) &= g\i\La^{G_+}(g)  
     = {\tsize\sum}_i X_i \otimes \Ad(g\i)\pr_+(\Ad(g)Y_i),\tag9\\ 
\de\la^{G_+}(e)X &= 0+{\tsize\sum}_i X_i \otimes \pr_+(\ad(X)Y_i),\\ 
\ga(\de\la^{G_+}(e)X,Y_k\otimes Y_l)  
&= {\tsize\sum}_i \ga(X_i,Y_k)\ga(\pr_+\ad(X)Y_i,Y_l)\\ 
&= \ga(\pr_+\ad(X)Y_k,Y_l) = \ga(\ad(X)Y_k,\pr_+^*Y_l)\\ 
&= \ga([X,Y_k],\pr_-Y_l) = \ga(X,[Y_k,Y_l]), 
\endalign$$ 
which we had to prove. 
Let us now investigate the Lie Poisson structure on $G$. From  
\nmb!{5.2}.\thetag2 we have 
$$\align 
\La_-(gu) &= {\tsize\sum}_i  
     \Bigl(guY_i\otimes guX_i-Y_igu\otimes X_igu\Bigr)\\ 
&= {\tsize\sum}_i g\Bigl(\Ad(u)Y_i\otimes \Ad(u)X_i 
     -\Ad(g\i)Y_i\otimes \Ad(g\i)X_i\Bigr)u \\ 
&= {\tsize\sum}_i g\Bigl(\Ad(u)Y_i\otimes \pr_-(\Ad(u)X_i)  
     + \Ad(u)Y_i\otimes \pr_+\Ad(u)X_i \\ 
&\quad-\pr_-(\Ad(g\i)Y_i)\otimes \Ad(g\i)X_i  
     -\pr_+(\Ad(g\i)Y_i)\otimes \Ad(g\i)X_i\Bigr)u. 
\endalign$$ 
In $L(\g,\g)$ we again have 
$$\align 
\Bigl({\tsize\sum}_i\Ad(u)Y_i\otimes &\pr_+\Ad(u)X_i\Bigr)\o \ga 
     = \pr_+\Ad(u)\pr_+\Ad(u\i) \\ 
&= \pr_+\Ad(u)(\Id-\pr_-)\Ad(u\i)  
     = \pr_+ - 0,\\ 
-\Bigl({\tsize\sum}_i \pr_-(\Ad(g\i)Y_i)\otimes &\Ad(g\i)X_i\Bigr)\o\ga 
     = -\Ad(g\i)\pr_+\Ad(g)\pr_-^* \\ 
&= -\Ad(g\i)\pr_+\Ad(g)\pr_+ = -\pr_+,\\ 
-\Bigl({\tsize\sum}_i \pr_+(\Ad(g\i)Y_i)\otimes &\Ad(g\i)X_i\Bigr)\o\ga  
     = -\Ad(g\i)\pr_+\Ad(g)\pr_+^* \\ 
&= -\Bigl({\tsize\sum}_i X_i\otimes\Ad(g\i)\pr_+(\Ad(g)Y_i)\Bigr)\o\ga.\\ 
\endalign$$ 
Thus we get 
$$\align 
\La_-(gu) &= {\tsize\sum}_i g\Bigl(\Ad(u)Y_i\otimes \pr_-(\Ad(u)X_i)  
     -X_i\otimes\Ad(g\i)\pr_+(\Ad(g)Y_i)\Bigr)u\\ 
&= g\Bigl({\tsize\sum}_i uY_i\otimes \pr_-(\Ad(u)X_i)u\Bigr)  
     -\Bigl({\tsize\sum}_i gX_i\otimes\pr_+(\Ad(g)Y_i)g\Bigr)u\\ 
&= g\La^{G_-}(u) - \La^{G_+}(g)u  
= T_{(g,u)}\ph(\La^{G_-}(u) - \La^{G_+}(g)),\\ 
\endalign$$ 
which proves \thetag6 and \thetag8.  
All remaining statements can be proved analogously, or are obvious. 
\qed\enddemo 
  
\proclaim{\nmb.{5.5}. Corollary} In the situation of \nmb!{5.1} we  
have: 
\roster 
\item The Poisson structure $\La^\ph_+$ on the direct product group 
       $G_+^{\text{op}}\x G_-$ is affine with  
$$\align 
(\La^\ph_+)_r(g,u) &= \La^{G_+}(g) + \La^{G_-}(u), \\ 
(\La^\ph_+)_l(g,u) &= \La^{G_+}(g) + \La^{G_-}(u) +  
     {\tsize\sum_i}Y_iu\wedge X_ig - {\tsize\sum_i}uY_i\wedge gX_i, 
\endalign$$ 
       where the vector fields $g\mapsto gX_i,X_ig$ are left and right  
       invariant with respect to the opposite  
       group structure on $G_+$. 
\item Moreover, the Lie-Poisson structure  
       $(\La^\ph_+)_r$ on $G_+^{\text{op}}\x G_-$ is the dual Lie  
       Poisson structure to $\La_-$ on $G$, i.e\., it defines the Lie  
       algebra structure on $\g$.   
\item The Poisson structure $\La^\ps_+$ on the direct product group 
       $G_-^{\text{op}}\x G_+$ is affine with 
$$\align 
(\La^\ps_+)_r(v,h) &= -\La^{G_+}(h) - \La^{G_-}(v), \\ 
(\La^\ps_+)_l(v,h) &= -\La^{G_+}(h) - \La^{G_-}(v) +  
     {\tsize\sum_i}Y_iv\wedge X_ih - {\tsize\sum_i}vY_i\wedge hX_i, 
\endalign$$ 
       where the vector fields $g\mapsto vY_i,Y_iv$ are left and  
       right invariant with respect to the 
       opposite group structure on $G_-$.  
\item $(G_+,-\La^{G_+})$ and $(G_-,\La^{G_-})$ are Lie-Poisson  
       subgroups of the Lie-Poisson group $(G,\La_-)$. 
\item The (local) projections from \nmb!{5.3} 
$$
p_l^+,p_r^+:(G,\La_-)\to (G_+,-\La^{G_+}),  \qquad
p_l^-,p_r^-:(G,\La_-)\to (G_-,\La^{G_-}),  
$$ 
     are Poisson mappings. 
\item The (local) projections from \nmb!{5.3} 
$$\alignat2 
&p_l^+:(G,\La_+)\to (G_+,\La^{G_+}),  
     &\qquad &p_r^+:(G,\La_+)\to (G_+,-\La^{G_+}), \\ 
&p_l^-:(G,\La_+)\to (G_-,\La^{G_-}),  
     &\qquad &p_r^-:(G,\La_+)\to (G_-,-\La^{G_-}) \\ 
\endalignat$$ 
     are Poisson mappings. 
\item The mapping $(G_+,\La^{G_+})\x (G,\La_+)\to (G,\La_+)$ given by  
       $(g,a)\mapsto ga$ is a left Poisson action of a Lie-Poisson  
       group. 
\item The mapping $(G,\La_+)\x (G_-,\La^{G_-})\to (G,\La_+)$ given by  
       $(a,u)\mapsto ga$ is a right Poisson action of a Lie-Poisson  
       group. 
\item The Lie-Poisson group dual to $(G,\La_-)$ is  
       $G_+\x G_-^{\text{op}}$ with the Lie-Poisson structure 
       $-(\La^\ps_+)_l$.
\endroster 
\endproclaim 
 
\demo{Proof} On the direct product group  
$G_+^{\text{op}}\x G_-$ the vector field $g\mapsto X_ig$ is  
{\it right} invariant, so the expressions in \therosteritem1 follows  
directly from from \nmb!{4.3} and the form \nmb!{5.4}.\therosteritem3  
of $\La^\ph_+$.   
The Poisson structure $(\La^\ph_+)_r$ is then visibly a Lie Poisson  
structure on $G_+^{\text{op}}\x G_-$, so $(\La^\ph_+)_r$ is affine.    
The proof of \therosteritem3 is similar. 
 
 For \therosteritem2 we consider, as explained in \nmb!{4.1} and in  
the proof of \nmb!{4.2}, see also the proof of \nmb!{5.4}.\thetag9: 
$$\align 
\la^\ph_{+,l}(g,u) &= \la^{G_+}(g) + \la^{G_-}(u)  
     + {\tsize\sum}_i \Ad(u\i)Y_i\wedge X_i  
     - {\tsize\sum}_i Y_i\wedge \Ad(g\i)X_i, \\ 
\de\la^\ph_{+,l}(e,e)(X,Y) &= \de\la^{G_+}(e)X + \de\la^{G_-}(e)Y  
     - {\tsize\sum}_i [Y,Y_i]_{\g_-}\wedge X_i  
     + {\tsize\sum}_i Y_i\wedge [X,X_i]_{\g_+} \\ 
&= b^{\g_-} + b^{\g_+}  
     - {\tsize\sum}_i [Y,Y_i]_{\g_-}\wedge X_i  
     + {\tsize\sum}_i Y_i\wedge [X,X_i]_{\g_+}, \\ 
\endalign$$ 
where $X\in\g_+$ and $Y\in\g_-$. If we take this into the inner  
product with elements $Y_k\otimes Y_l$, $Y_k\otimes X_l$, etc., use  
\nmb!{5.4}.\therosteritem9 and proceed as there, the result follows. 

\therosteritem5 -- \therosteritem8 follow from the formulae for  
$\La_+$ and $\La_-$ in the `coordinates' $(g,u)$ and $(v,h)$, and  
from the fact that $\La^{G_+}$ and $\La^{G_-}$ are multiplicative. 
 
\therosteritem9 is analogous to \therosteritem2.
\qed\enddemo 
 
\subhead\nmb.{5.6} \endsubhead 
Let us note finally that the decompositions \nmb!{5.4}~\thetag3 and  
\thetag4 of the Poisson structure $\La_+$ on $G\cong G_+\x G_-$ are  
surprisingly rigid. 
 
\proclaim{Theorem} Suppose that a Poisson structure $\La$ on a  
manifold $H\x K$ which is a product of two Lie groups of equal  
dimension admits a decomposition 
$$ 
\La(h,k) = \La^{H}(h) + \La^{K}(k)  
     + {\tsize\sum}_i Y_i^r(k)\wedge  X_i^l(h)  
     \in{\tsize\bigwedge}^2T_{(h,k)}(H\x K), 
$$ 
where $\La^H$ and $\La^K$ are tensor fields on $H$ and $K$,  
respectively, and where $X_i^l$ are the left invariant vector fields  
and $Y_i^r$ the right invariant vector fields on $H$ and $K$, with  
respect to bases $X_i$ of $\frak h$ and $Y_i$ of $\frak k$. 
 
Then $\La^H$ and $\La^K$ are affine Poisson structures on $H$ and  
$K$, respectively, and $(H,\La^H)$, $(K,\La^K)$ is a dual pair of  
Lie-Poisson groups and $\La$ represents the `symplectic' Poisson  
tensor on the corresponding group double. 
\endproclaim 
 
\demo{Proof} 
The vanishing Schouten bracket $[\La,\La]$ yields 
$$\alignat2 
0 &= [\La^H,\La^H]          &\qquad& \in \Ga({\tsize\bigwedge}^3 TH) \\ 
&\quad + [\La^K,\La^K]      &\qquad& \in \Ga({\tsize\bigwedge}^3 TK) \\ 
&\quad + 2{\tsize\sum}_iY_i^r\wedge [X_i^l,\La^H]  
     - {\tsize\sum}_{ij} [Y_i,Y_j]^r\wedge X_i^l\wedge X_j^l   
     &\qquad& \in \X(K)\otimes\Ga({\tsize\bigwedge}^2 TH) \\ 
&\quad - 2{\tsize\sum}_i[Y_i^r,\La^K]\wedge X_i^l  
     + {\tsize\sum}_{ij} Y_i^r\wedge Y_j^r\wedge [X_i, X_j]^l 
      &\qquad& \in \Ga({\tsize\bigwedge}^2 TK)\otimes \X(H). \\ 
\endalignat$$ 
Each of the lines vanishes by itself: The first two lines then say  
that $\La^H$ and $\La^K$ are Poisson tensors on $H$ and $K$,  
respectively.  
Using the structure constants $c^{ij}_m$ of $\frak h$ with respect to  
the basis $X_i$, and $d^{ij}_m$ of $\frak k$ with respect to $Y_i$,  
the last two lines can be rewritten as  
$$\align 
{\tsize\sum}_mY_m^r\wedge [X_m^l,\La^H]  
     &=\tfrac{1}2 {\tsize\sum}_{ijm}  
     d^{ij}_mY_m^r\wedge X_i^l\wedge X_j^l, \\ 
{\tsize\sum}_m[Y_m^r,\La^K]\wedge X_m^l  
     &=\tfrac12 {\tsize\sum}_{ijm} Y_i^r\wedge Y_j^r\wedge c^{ij}_m 
     X_m^l, 
\endalign$$ 
or by  
$$\align 
[X_m^l,\La^H] &= \L_{X_m^l}\La^H  
     =\tfrac{1}2 {\tsize\sum}_{ij}  
     d^{ij}_m X_i^l\wedge X_j^l \in L({\tsize\bigwedge}^2\frak h),\\ 
[Y_m^r,\La^K] &= \L_{Y_m^r}\La^K 
     =\tfrac12 {\tsize\sum}_{ij} c^{ij}_m Y_i^r\wedge Y_j^r  
     \in R({\tsize\bigwedge}^2\frak k). 
\endalign$$ 
These are just conditions \therosteritem3 and \therosteritem4 of  
\nmb!{4.1} without the further assumption that $\La^H(e)=0$ or  
$\La^K(e)=0$, so we can conclude from there that $\La^H$ and $\La^K$  
are affine Poisson structures, respectively.  
For their associated Lie-Poisson structures  
$$\align 
(\La^H)_r(h) &= \La^H(h)-\La^H(e)h,\\ 
(\La^K)_l(k) &= \La^K(h)-k\La^K(e) 
\endalign$$ 
we get  
$$\align 
\L_{X_m^l}(\La^H)_r &= \L_{X_m^l}\La^H
     =\tfrac{1}2 {\tsize\sum}_{ij} d^{ij}_m X_i^l\wedge X_j^l, \\ 
\L_{Y_m^r}(\La^K)_l &= \L_{Y_m^r}\La^K 
     =\tfrac12 {\tsize\sum}_{ij} c^{ij}_m Y_i^r\wedge Y_j^r,  
\endalign$$ 
so that the Lie Poisson structure $(\La^H)_r$ corresponds to the  
cobracket 
$$  
b_{\frak h}':\frak h\to {\tsize\bigwedge}^2\frak h,\quad 
     b_{\frak h}'(X_m)  
     =\tfrac{1}2 {\tsize\sum}_{ij} d^{ij}_m X_i\wedge X_j,  
$$ 
and the Lie Poisson bracket $(\La^K)_l$ corresponds to the cobracket  
$$  
b_{\frak k}':\frak k\to {\tsize\bigwedge}^2\frak k,\quad 
     b_{\frak k}'(Y_m)  
     =\tfrac{1}2 {\tsize\sum}_{ij} c^{ij}_m Y_i\wedge Y_j. 
$$ 
Hence $b_{\frak h}'$ is dual to the Lie bracket on $\frak k$, and  
$b_{\frak k}'$ is dual to the Lie bracket on $\frak h$, with respect  
to the pairing $\ga(X_i,Y_j)=\de_{ij}$.  
\qed\enddemo 
 
\head\totoc\nmb0{6}. Dressing Actions and Symplectic Leaves \endhead 
 
\subhead\nmb.{6.1}. Lie algebroids\endsubhead 
On every Poisson manifold $(M,\La)$ the Poisson tensor defines the  
mapping $T^*M\ni \al\mapsto \al^\sharp:=\bar\imath_\al\La\in TM$, and  
a Lie bracket on the space of 1-forms defined by  
$$ 
\{\al,\be\} := i_{\al^\sharp}d\be - i_{\be^\sharp}d\al +  
di_\La(\al\wedge \be). 
\tag1$$ 
The mapping $(\quad)^\sharp:\Om^1(M)\to \X(M)$ is then a homomorphism  
of Lie algebras,  
$$ 
\{\al,\be\}^{\sharp} = [\al^\sharp,\be^\sharp]; 
\tag2$$ 
this is also expressed by saying that $\La$ turns $T^*M$ into a Lie  
algebroid with anchor mapping $(\quad)^\sharp$.  
 
\subhead\nmb.{6.2}. The dressing action  \endsubhead 
Affine Poisson structures on a Lie group $G$ may be characterized by  
the property that the left invariant 1-forms (or equivalently the  
right invariant ones) are closed with respect to the bracket  
\nmb!{6.1}.\thetag1.   
 
Consequently, for an affine Poisson structure $\La$ on $G$ the  
mappings 
$$\alignat2 
\la:\g^*&\to \X(G),&\quad \la(X)(a) &:= -(aX)^\sharp,\\ 
\rh:\g^*&\to \X(G),&\quad \rh(X)(a) &:= (Xa)^\sharp\\ 
\endalignat$$ 
are an anti homomorphism and homomorphism of the Lie algebras  
$\g^*_r$ and $\g^*_l$, respectively, where $\g^*_r$ is the dual space  
$\g^*$ with the Lie bracket corresponding to $\La_r$, and where  
$\g^*_l$ corresponds to $\La_l$.  
The fields $\la(X)$ are called  
\idx{\it left dressing vector fields} on $G$, and the $\rh(X)$ are  
called \idx{\it right dressing vector fields}. They may be considered  
as infinitesimal actions of the corresponding dual groups. We have  
seen such actions already in  
\nmb!{5.5}.\therosteritem4 and \therosteritem5. 
If we can integrate this infinitesimal action to a global one, called  
the dressing action (if the  
dressing fields are complete), the affine Poisson group $(G,\La)$  
will be called \idx{\it complete}.  
 
In any case, the left (or right) dressing vector fields generate the  
characteristic distribution of $\La$, whose leaves are precisely the  
symplectic leaves of the Poisson structure $\La$.  
 
One believes that dressing actions describe `hidden symmetries' of  
physical systems.  
 
\proclaim{\nmb.{6.3}. Theorem} Let $G$ be a Lie group with a metrical  
Lie algebra $(\g,\ga)$ which admits a Manin decomposition  
$\g=\g_+\oplus \g_-$. In the setting of \nmb!{5.1}, the dressing  
vector fields for the affine Poisson structures $\La_+$ and $\La_-$  
on $G$ are the following: 
$$\alignat2 
\la_+(X_i)(a) &= -\pr_+(\Ad(a)X_i)a, &\qquad  
     \rh_+(X_i)(a) &= a\pr_+(\Ad(a)\i X_i),\tag1\\  
\la_+(Y_i)(a) &= \pr_-(\Ad(a)Y_i)a, &\qquad  
     \rh_+(Y_i)(a) &= -a\pr_-(\Ad(a)\i Y_i).\\  
\la_-(X_i)(a) &= \pr_-(\Ad(a)X_i)a, &\qquad  
     \rh_-(X_i)(a) &= -a\pr_-(\Ad(a)\i X_i),\tag2\\  
\la_-(Y_i)(a) &= -\pr_+(\Ad(a)Y_i)a, &\qquad  
     \rh_-(Y_i)(a) &= a\pr_+(\Ad(a)\i Y_i).\\  
\endalignat$$ 
\endproclaim 
 
\demo{Proof} For instance, by \nmb!{5.2}, 
$$\align 
\rh_+(X_i)(a) &= \bar\imath(\ga(X_ia))\La_+(a)  
     = \bar\imath(\ga(X_ia)) 
     {\tsize\sum}_j(Y_ja\otimes X_ja-aX_j\otimes aY_j)\\ 
&= {\tsize\sum}_j(\ga(X_ia,Y_ja)X_ja-\ga(X_ia,aX_j)aY_j)\\ 
&= X_ia-a\pr_-(\Ad(a\i)X_i 
     = a(\Ad(a\i)X_i - \pr_-(\Ad(a\i)X_i)) \\ 
&= a(\pr_+(\Ad(a\i)X_i)). \qed 
\endalign$$ 
\enddemo 
 
\proclaim{\nmb.{6.4}. Corollary} The Poisson tensors $\La_\pm$ may be  
written in the following alternative form: 
$$\align 
\La_+(a) &= {\tsize\sum}_i \Bigl(aY_i\oplus \pr_+(\Ad(a)X_i)a  
     -aX_i\otimes\pr_-(\Ad(a)Y_i)a\Bigr) \tag1\\ 
&= {\tsize\sum}_i \Bigl(Y_ia\oplus a\pr_+(\Ad(a\i)X_i)  
     -X_ia\otimes a\pr_-(\Ad(a\i)Y_i)\Bigr). \\ 
\La_-(a) &= -{\tsize\sum}_i \Bigl(aX_i\oplus \pr_+(\Ad(a)Y_i)a  
     -aY_i\otimes\pr_-(\Ad(a)X_i)a\Bigr) \tag2\\ 
&= -{\tsize\sum}_i \Bigl(X_ia\otimes a\pr_-(\Ad(a\i)Y_i) 
     -Y_ia\oplus a\pr_+(\Ad(a\i)X_i)\Bigr). \\ 
\endalign$$  
\endproclaim 
 
\demo{Proof} 
 From the definition of $(\quad)^\sharp:T^*G\to TG$ we have 
$$\La_+ =  
{\tsize\sum}_i\Bigl(-aY_i\otimes\la_+(X_i)-aX_i\otimes\la_+(Y_i)\Bigr), 
\text{ etc. } \qed 
$$ 
\enddemo 
 
\proclaim{\nmb.{6.5}. Corollary} \cit!{1} 
\roster 
\item The characteristic distributions $S_\pm$ of the Poisson  
       structures $\La_\pm$ may be described as follows. 
$$\align 
S_+(a) &= a(\pr_+(\Ad(a\i)\g_+)+\pr_-\Ad(a\i)\g_-)\\ 
&= (\pr_+(\Ad(a)\g_+)+\pr_-\Ad(a)\g_-)a, \\ 
S_-(a) &= a(\pr_-(\Ad(a\i)\g_+)+\pr_+\Ad(a\i)\g_-)\\ 
&= (\pr_-(\Ad(a)\g_+)+\pr_+\Ad(a)\g_-)a. \\ 
\endalign$$  
       In particular, $S_+(a)+ S_-(a)=T_aG$. 
\item The symplectic leaves of $S_+$ are the connected components of  
       the intersections of orbits $G_+\,a\,G_-\cap G_-\,a\,G_+$, and  
       the symplectic leaves of $S_-$ are the connected components of  
       the intersections of orbits $G_-\,a\,G_-\cap G_+\,a\,G_+$, for  
       $a\in G$. 
\item The Poisson structure $\La_+$ is non-degenerate precisely on  
       the set $G_+G_-\cap G_-G_+$; so it is globally non-degenerate  
       if and only if $G_+G_-=G$. In particular, if $(G,\La_+)$ is  
       complete (\nmb!{6.2}) then $\La_+$ is non-degenerate. 
\endroster 
\endproclaim 
 
\demo{Proof} 
\therosteritem1 follows directly from theorem \nmb!{6.3} since the  
dressing vector fields generate the characteristic distribution.  
To prove \therosteritem2 observe that the tangent space to the  
intersection of orbits $G_+\,a\,G_-\cap G_-\,a\,G_+$ at $a\in G$ is  
$$\align 
(\g_+a + a\g_-)\cap(\g_-a + a\g_+)  
     &= a((\Ad(a\i)\g_+ + \g_-)\cap(\Ad(a\i)\g_- + \g_+)) \\ 
     &= a(\pr_+(\Ad(a\i)\g_+) + \pr_-(\Ad(a\i)\g_-)) = S_+(a), 
\endalign$$ 
so that the connected components of $G_+\,a\,G_-\cap G_-\,a\,G_+$ are  
integral submanifolds of $S_+$. For $S_-$ the proof is similar. 
 
The intersection $G_+G_-\cap G_-G_+$ is an open and dense subset of  
$G$ consisting, by \therosteritem2, of points where $\La_+$ is  
non-degenerate. If the orbit $G_+aG_-$ meets $G_+G_-\cap G_-G_+$ then  
it is contained in $G_+G_-$, so $G_+G_-\cap G_-G_+$ consists of all  
points where $\La_+$ is non-degenerate.  
\qed\enddemo 
 
\subhead\nmb.{6.6} \endsubhead 
On $M:= G_+G_-\cap G_-G_+$ the Poisson structure $\La_+$ is  
symplectic, so let us describe the associated symplectic form  
$\om=(\La_+)\i$ in terms of the coordinates $(g,u)$ and $(v,h)$  
introduced in \nmb!{5.3}. We will start by describing the dressing  
vector fields on the groups $(G_+\x G_-, \La^\ph_+)$ and $(G_-\x G_+,  
\La^\ps_+)$. In order to avoid problems of always having to tell  
which multiplication is opposite, and to use a notation which differs  
from that used in theorem \nmb!{6.3} we will write $(uX)^\sharp$ for  
the dressing vector field corresponding to the left invariant 1-form  
on $G_+\x G^*$ represented by $\et(g,u)=uX$ in the obvious way: 
$$ 
\ga(uX,gX_i+uY_j) = \ga(uX,uY_j) = \ga(X,Y_j),\quad\text{ etc.} 
$$  
After easy calculations we get from \nmb!{5.4}.\thetag3 and \thetag4: 
 
\proclaim{Theorem} In the situations above, the dressing vector  
fields are given by: 
$$\align 
\text{ On }(G_+\x G_-&,\La^\ph_+(g,u)):\quad\tag1\\ 
(X_iu)^\sharp &= gX_i - u\pr_-(\Ad(u\i)X_i), \\ 
(Y_ig)^\sharp &= -g\pr_+(\Ad(g\i)Y_i) - \pr_-(\Ad(g\i)Y_i)u, \\ 
(uX_i)^\sharp &= \pr_-(\Ad(u)X_i)u + g\pr_+(\Ad(u)X_i)u, \\ 
(gY_i)^\sharp &= \pr_+(\Ad(g)Y_i)g - Y_iu. \\ 
\text{ On }(G_-\x G_+&,\La^\ps_+(v,h)):\quad\tag2\\ 
(vX_i)^\sharp &= X_ih - \pr_-(\Ad(v)X_i)v, \\ 
(hY_i)^\sharp &= -\pr_+(\Ad(h)Y_i)h - v\pr_-(\Ad(h)Y_i), \\ 
(X_iv)^\sharp &= v\pr_-(\Ad(v\i)X_i) + \pr_+(\Ad(v\i)X_i)h, \\ 
(Y_ih)^\sharp &= h\pr_+(\Ad(h\i)Y_i) - vY_i. \\ 
\endalign$$ 
Denote now $(X_iu)^\wedge = \ph_*(X_u)^\sharp \in \X(G_+G_-)$, etc., and  
$(vX_i)^\wedge = \ps_*(X_iu)^\sharp\in \X(G_-G_+)$, etc., and call them  
the \idx{\rm undressing vector fields}. They  
are given at the point $a=gu=vh\in M= G_+G_-\cap G_-G_+ \subset G$ by  
$$\alignat2 
(X_iu)^\wedge &= a\pr_+(\Ad(u\i)X_i), &\qquad 
     (Y_ig)^\wedge &= -Y_ia, \tag3\\ 
(uX_i)^\wedge &= aX_i, &\qquad 
     (gY_i)^\wedge &= -\pr_-(\Ad(g)Y_i)a, \\ 
(vX_i)^\wedge &= \pr_+(\Ad(v)X_i)a &\qquad 
     (hY_i)^\wedge &= -aY_i, \\ 
(X_iv)^\wedge &= X_ia, &\qquad 
     (Y_ih)^\wedge &= -a\pr_-(\Ad(h\i)Y_i).\\ 
\endalignat$$ 
\endproclaim 
 
\demo{Proof} 
We only prove \thetag3, and only one example: 
$$\align 
(X_iu)^\wedge &= \ph_*(X_iu)^\sharp = \ph_*(gX_i-u\pr_-(\Ad(u\i)X_i))\\ 
&=  gX_iu-gu\pr_-(\Ad(u\i)X_i) = gu(\Ad(u\i)X_i -\pr_-(\Ad(u\i)X_i))\\ 
&= a\pr_+(\Ad(u\i)X_i). \qed 
\endalign$$ 
\enddemo 
 
\proclaim{\nmb.{6.7}. Corollary} At points  
$a=gu=vh\in M= G_+G_-\cap G_-G_+ \subset G$ the affine Poisson  
structure is given by 
$$\align 
\La_+(a) &= {\tsize\sum}_i((uX_i)^\wedge \otimes (hY_i)^\wedge  
     + (X_iv)^\wedge \otimes (Y_ig)^\wedge) \tag1\\ 
&= {\tsize\sum}_i((X_iu)^\wedge \otimes (gY_i)^\wedge  
     - (Y_ih)^\wedge \otimes (vX_i)^\wedge). \\ 
\endalign$$ 
The associated symplectic structure $\om$ may be written as 
$$\align 
\om_a &= {\tsize\sum}_i((uX_i)\otimes (hY_i) 
     + (X_iv)\otimes (Y_ig))\tag2\\ 
&= {\tsize\sum}_i((X_iu)\otimes (gY_i) 
     - (Y_ih)\otimes (vX_i)^\wedge), \\ 
\endalign$$ 
where we identify the 1-forms $uX_i$, etc., on $G_+\x G_-$ and the  
1-forms $hY_i$, etc., on $G_-\x G_+$ with 1-forms on $M$ via the  
diffeomorphisms $\ph$ and $\ps$. Formally correct we should write  
$(\ph\i)^*(uX_i)$, etc. 
\endproclaim 
 
\demo{Proof} 
The form \thetag1 of $\La_+(a)$ can be checked by easy calculations.  
But \thetag1 shows that we can construct $\La_+(a)$ from  
$(uX_i)^\wedge = \bar\imath(uX_i)\La_+$, etc., thus we can construct  
$\om_a=\La_+(a)\i$ in the same way from the corresponding 1-forms  
$uX_i$.  
\qed\enddemo 
 
\subhead\nmb.{6.8}. Remark \endsubhead 
We can write \nmb!{6.7}.\therosteritem2 in a more `coordinate free'  
form: 
$$\align 
\om &= \ga(\mu^\ph_{G_-}\;\overset\otimes \to ,\; \mu^\ps_{G_+}) 
    + \ga(\th^\ps_{G_-}\;\overset\otimes \to ,\; \mu^\ph_{G_+})\tag1\\ 
&= \frac12(\ga(\th^\ph_{G_-}\;\overset\wedge \to ,\; \mu^\ph_{G_+}) 
    + \ga(\mu^\ps_{G_-}\;\overset\wedge \to ,\; \th^\ps_{G_+})), 
\endalign$$ 
where $\mu^\ph_{G_-}= (uX_i)\otimes Y_i$ is the left Maurer Cartan  
form on $G_-$ pushed via $\ph$ to $M= G_+G_-\cap G_-G_+ \subset G$,  
and where $\th^\ps_{G_-}= (X_iv)\otimes Y_i$ is the right Maurer Cartan  
form on $G_-$ pushed via $\ps$ to $M$, etc. This expression \thetag1  
should be compared with the corresponding formula in \cit!{1}, or  
with formula 2.3.(3) in \cit!{2} for the case of a cotangent bundle  
$T^*G_+$. So \nmb!{6.7} is a generalization of these results in  
\cit!{2} to the case of a double group. 
 
\subhead\nmb.{6.9} \endsubhead 
Recall now from \nmb!{5.3} the   
projections $p^+_l,p^+_r:G\supset U\to G_+$ and  
$p^-_l,p^-_r:G\supset U\to G_-$ which we get from inverting $\ph$ and  
$\ps$, respectively. For $a\in G$ and for $b$ near $e$ in $G$ we then  
define  
$$\alignat2 
\la^+_b(a) &:= p_r^+(ab\i a\i)a, &\qquad  
     \la^-_b(a) &:= p_r^-(ab\i a\i)a, \tag1\\ 
\rh^+_b(a) &:= ap_l^+(ab\i a\i), &\qquad  
     \rh^-_b(a) &:= ap_l^-(ab\i a\i). \\ 
\endalignat$$ 
 
\proclaim{Theorem} 
The mappings $\la^+$ and $\la^-$ define left (local) actions of $G$  
on $G$, and $\rh^+$ and $\rh^-$ define right (local actions), i.e., 
$$\alignat2 
\la^+_b(\la^+_{b'}(a)) &= \la^+_{bb'}(a), &\qquad  
     \la^-_b(\la^-_{b'}(a)) &= \la^-_{bb'}(a), \tag2\\ 
\rh^+_b(\rh^+_{b'}(a)) &= \rh^+_{b'b}(a), &\qquad  
     \rh^-_b(\rh^-_{b'}(a)) &= \rh^-_{b'b}(a).\\ 
\endalignat$$ 
The subgroup $G_+$ is invariant under $\la^+$ and $\rh^+$ while  
$G_-$ is invariant under $\la^-$ and $\rh^-$. 
 
Moreover, the pairs $\la^+$, $\la^-$, and $\rh^+$, $\rh^-$ commute: 
$$\align 
\la^+_b(\la^-_{b'}(a)) &= \la^-_{b'}(\la^+_{b}(a)) =  
     p^+_r(a(b')\i ba\i)ab\i, \tag3 \\ 
\rh^+_b(\rh^-_{b'}(a)) &= \rh^-_{b'}(\rh^+_{b}(a)) =  
     (b')\i a p^+_r(a(b')\i ba\i).\\ 
\endalign$$  
\endproclaim 
  
\demo{Proof} \thetag2. 
Assume $a(b')\i a\i = vh$ for $v\in G_-$ and $h\in G_+$ as usual, so  
that we have $\la^+_{b'}(a)= ha = v\i a(b')\i$. Then 
$$\align 
\la^+_b(\la^+_{b'}(a)) &= p^+_r(hab\i a\i h\i)ha =  p^+_r(hab\i a\i)a  
     = p^+_r(v\i a(b')\i b\i a\i)a \\ 
&= p^+_r(a(bb')\i a\i)a = \la^+_{bb'}(a). 
\endalign$$ 
For the other actions, the proofs are similar.  
 
\thetag3. We shall prove only the first part.  
Put $a(b')\i a\i = gu$ for $g\in G_+$ and $u\in G_-$, so that  
$\la^-_{b'}(a) = ua = g\i a (b')\i$. Then 
$$\align 
\la^+_b(\la^-_{b'}(a)) &= \la^+_b(ua) = p^+_r(uab\i (ua)\i)ua  
     = p^+_r(ab\i b' a\i g)g\i a(b')\i\\  
&= p^+_r(ab\i b' a\i)a(b')\i  
     = p^-_r(ab\i b' a\i)\i (ab\i b' a\i) a(b')\i\\  
&= p^-_r(a(b')\i b a\i) ab\i.  
\endalign$$ 
On the other hand, put $ab\i a\i = vh$, so that  
$\la^+_b(a)=ha= v\i a b\i$. Then 
$$\align 
\la^-_{b'}(\la^+_{b}(a)) &= p^-_r(ha(b')\i (ha)\i)ha  
     = p^-_r(a(b')\i b a\i v)v\i ab\i \\ 
&= p^-_r(a(b')\i b a\i) ab\i. \qed 
\endalign$$ 
\enddemo 
 
\proclaim{\nmb.{6.10}. Theorem} 
The infinitesimal actions for $\la^+$, $\la^-$, $\rh^+$, and $\rh^-$  
are the following, where $A,B\in\g=\g_+\oplus\g_-$ and $a\in G$: 
$$\alignat2 
\la^+_B(a) &:= -\pr^+(\Ad(a)B)a, &\qquad  
     \la^-_B(a) &:= -\pr^-(\Ad(a)B)a, \tag1\\ 
\rh^+_B(a) &:= -a\pr^+(\Ad(a\i)B), &\qquad  
     \rh^-_B(a) &:= -a\pr^-(\Ad(a\i)B). 
\endalignat$$ 
Furthermore, as usual for left and right actions, for $B, B'\in \g$ we  
have 
$$\alignat2 
[\la^+_B,\la^+_{B'}] &= -\la^+_{[B,B']}, &\qquad  
     [\la^-_B,\la^-_{B'}] &= -\la^-_{[B,B']},\tag2\\ 
[\rh^+_B,\rh^+_{B'}] &= \rh^+_{[B,B']}, &\qquad  
     [\rh^-_B,\rh^-_{B'}] &= \rh^-_{[B,B']},\\ 
[\la^+_B,\la^-_{B'}] &= 0, &\qquad  
     [\rh^+_B,\rh^-_{B'}] &= 0.\\ 
\endalignat$$ 
Moreover,  
$$\alignat2 
\la_+(X_i) &= \la^+_{X_i}, &\quad \la_+(Y_i) &= -\la^-_{Y_i},\tag3\\ 
\rh_+(X_i) &= -\rh^+_{X_i}, &\quad \rh_+(Y_i) &= -\rh^-_{Y_i}, \\ 
\la_-(X_i) &= -\la^-_{X_i}, &\quad \la_-(Y_i) &= \la^+_{Y_i}, \\ 
\rh_-(X_i) &= -\rh^-_{X_i}, &\quad \rh_-(Y_i) &= \rh^+_{Y_i},  
\endalignat$$ 
so that we can reconstruct the dressing actions from  
$\la^+$, $\la^-$, $\rh^+$, and $\rh^-$. For example, the (local) left  
dressing action for $\La_+$ is given by  
$$\gather 
G_+ \x (G_-)^{\text{op}} \x G \to G, \\ 
(g,u).a = \la^+_g\la^-_{u\i}(a) = p^-_r(auga\i)ag\i. 
\endgather$$  
The (local) left dressing action for $\La_-$ is given by  
$$\gather 
G_+ \x (G_-)^{\text{op}} \x G \to G, \\ 
(g,u).a = \la^-_g\la^+_{u\i}(a) = p^-_r(ag\i u\i a\i)au. \qed 
\endgather$$  
\endproclaim 
 
\subhead\nmb.{6.11}. Remark \endsubhead 
The dressing actions of $G_+$ on $G_-$, and of $G_-$ on $G_+$ can  
also be reconstructed from this scheme. For $(g,u)\in G_+\x G_-$  
they are given by restricting the (local) actions  
$\la^+$, $\la^-$, $\rh^+$, and $\rh^-$ of $G$ on $G$ appropriately  
(see \nmb!{6.9}): 
$$\alignat2 
\la^+_u(g) &:= p_r^+(gu\i g\i)g, &\qquad  
     \la^-_g(u) &:= p_r^-(ug\i u\i)u,\\ 
\rh^+_u(g) &:= gp_l^+(gu\i g\i), &\qquad  
     \rh^-_g(u) &:= up_l^-(ug\i u\i). \\ 
\endalignat$$ 
Note that in these formulae one should replace, say,  
$p_r^+(gu\i g\i)g=p^+_r(gu\i)$ only if the action is complete, or  
only for $g$ and $u$ near $e$, since the left hand side is defined  
for all $g$ and for $u$ near $e$, whereas the right hand side needs  
both $g$ and $u$ near $e$.  
 
\proclaim{\nmb.{6.12}. Corollary} 
The dressing actions of $G_+$ on $G_-$, and of $G_-$ on $G_+$ are  
(local) Poisson actions. 
\endproclaim 
 
\demo{Proof} 
We prove it only, say, for the left dressing action of $G^-$ on  
$G_+$. At least locally this is given by 
$$ 
\la^+:G_-\x G_+ \to G_+,\qquad \la^+(u,g) = \la^+_u(g) = p^+_r(gu\i), 
$$ 
Due to theorem \nmb!{5.4}.\thetag5 the mapping  
$$ 
\tilde\ph:G_+\x G_-\ni (g,u)\mapsto gu\i\in G 
$$ 
is a Poisson mapping  
$(G_+\x G_-,\La^{G_+}\x \La^{G_-})\to (G,\La_-)$. 
By corollary \nmb!{5.5} the (local) projection  
$p^+_r:(G,\La_-)\to (G^+, -\La^{G_+})$ is a Poisson mapping, so the  
composition $\la^+ = p^+_r\o \tilde\ph$ is also Poisson. 
\qed\enddemo 
 
\head\totoc\nmb0{7}. Examples \endhead 
 
\subhead\nmb.{7.1}. Example \endsubhead 
Let us assume that in the Manin decomposition $\g=\g_+\oplus\g_-$ the  
subalgebra $\g_-$ is commutative then the simply connected Lie group  
$G$ is isomorphic to the cotangent bundle  
$T^*G_+\cong G_+\ltimes \g_-$, the semidirect product of $G_+$ and  
the dual Lie algebra, which is complete with respect to the dressing  
actions. $\ph: G_+\x \g_-\to T^*G_+$ is the left trivialization,  
$\ps$ is the right trivialization. This situation was described in 
detail in our earlier paper \cit!{2}.
 
\subhead\nmb.{7.2}. Example \endsubhead 
We consider $\g_+ = \frak s\frak u(2)$ with the standard matrix basis 
$$ 
e_1=\frac12\pmatrix i & 0 \\ 0 & -i \endpmatrix, \quad 
e_2=\frac12\pmatrix 0 & 1 \\ -1 & 0 \endpmatrix, \quad 
e_3=\frac12\pmatrix 0 & i \\ i & 0 \endpmatrix, 
$$ 
satisfying $[e_1,e_2]=e_3$, $[e_2,e_3]=e_1$, and $[e_3,e_1]=e_2$. 
The following commutation rules  
$[e^*_1,e^*_2]=e^*_2$, $[e^*_1,e^*_3]=e^*_3$, and $[e^*_2,e^*_3]=0$ 
for the dual basis in $\g_-=\g_+^*$  
make $\g=\g_+\oplus\g_-$ into a Lie bialgebra which is isomorphic to  
$\frak s\frak l(2,\Bbb C)$ as 6 dimensional real algebra with  
$\g_-=\frak s\frak b(2,\Bbb C)$, where the  
elements of the dual basis are given by  
$$ 
e_1^*=\frac12\pmatrix 1 & 0 \\ 0 & -1 \endpmatrix, \quad 
e_2^*=\pmatrix 0 & -i \\ 0 & 0 \endpmatrix, \quad 
e_3^*=\pmatrix 0 & 1 \\ 0 & 0 \endpmatrix. 
$$ 
The invariant symmetric pairing can be recognized as  
$$ 
\ga(A,B) = 2 \operatorname{Im}\tr(AB). 
$$ 
We consider now the double Lie group $G=SL(2,\Bbb C)$ with  
$G_+=SU(2)$ and $G_-=SB(2,\Bbb C)$. We will write the elements as  
follows: 
$$\align 
&G=SL(2,\Bbb C)\ni a = \pmatrix z_1 & z_2 \\ z_3 & z_4 \endpmatrix,  
     \quad \text{ where } z_i\in \Bbb C, \quad z_1z_4-z_2z_3 = 1, \\ 
&G_+=SU(2)\ni g = \pmatrix \al & -\bar\nu \\ \nu & \bar \al \endpmatrix,  
     \quad \text{ where }\al, \nu \in \Bbb C, \quad |\al|^2+|\nu|^2=1, \\ 
&G_-=SB(2,\Bbb C)\ni u =\pmatrix t& \ga \\ 0 & t\i\endpmatrix, 
     \quad \text{ where }t>0,  \ga\in \Bbb C. 
\endalign$$ 
We define $\La_+(a)=\tfrac12(ra+ar)$ with  
$r=\sum_ie^*_i\wedge e_i$ on $SL(2,\Bbb C)$ as explained in  
\nmb!{5.1}. We then extend it onto the whole space $GL(2,\Bbb C)$ of  
all invertible matrices by admitting $a\in GL(2,\Bbb C)$. Since the  
left and right invariant vector fields on  
$\frak g\frak l(2,\Bbb C)\cong \Bbb C^4\cong \Bbb R^8$ satisfy the  
same commutation rules as their restrictions to $SL(2,\Bbb C)$, we  
will get a Poisson structure. Of course it is tangent to  
$SL(2,\Bbb C)$, so that, if we consider the Poisson brackets between  
all matrix elements $z_i$ and $\bar z_i$, the functions  
$\det = z_1z_4 - z_2z_3$  and  
$\overline{\det} = \overline{z_1z_4} - \overline{z_2z_3}$ will be  
Casimirs for the bracket. Thus we get a Poisson structure on  
$GL(2,\Bbb C)$ whose restriction to $SL(2,\Bbb C)$ is exactly  
$\La_+$. We calculated the following Poisson brackets, which were  
also obtained independently by \cit!{36}. 
$$\alignat2 
\{z_1,z_2\} &= -\tfrac12 iz_1z_2 &\qquad \{z_2,z_3\} &= iz_1z_4 \\ 
\{z_1,z_3\} &= \tfrac12 iz_1z_3 &\qquad \{z_2,z_4\} &=\tfrac12iz_2z_4\\ 
\{z_1,z_4\} &= 0 &\qquad               \{z_3,z_4\} &=-\tfrac12iz_3z_4\\ 
\{z_1,\bar z_1\} &= -\tfrac12 i|z_1|^2 -i|z_3|^2 &\qquad  
     \{z_2,\bar z_2\} &=-\tfrac12 i|z_2|^2 -i|z_1|^2-i|z_4|^2  \\ 
\{z_3,\bar z_3\} &= -\tfrac12 i|z_3|^2 &\qquad  
     \{z_4,\bar z_4\} &= -\tfrac12i|z_4|^2 -i|z_3|^2\\ 
\{z_1,\bar z_2\} &= -iz_3\bar z_4 &\qquad  
     \{z_2,\bar z_3\} &=\tfrac12 iz_2\bar z_3 \\ 
\{z_1,\bar z_3\} &= 0 &\qquad \{z_2,\bar z_4\} &= -iz_1\bar z_3\\ 
\{z_1,\bar z_4\} &= \tfrac12 iz_1\bar z_4 &\qquad \{z_3,\bar z_4\} &= 0 
\endalignat$$ 
The lacking commutators may be obtained from this list if we remember  
that the Poisson bracket is real, e.g.,  
$\{\bar z_i,\bar z_j\} = \overline{\{z_i,z_j\}}$. 
One can then check that indeed $\det$ and $\overline{\det}$ are  
Casimir functions, and that $z_1\leftrightarrow z_4$, 
$z_2\mapsto -z_2$, and  
$z_3\mapsto -z_3$ 
defines a symmetry of the bracket associated to the inverse  
$a\mapsto a\i$ in $SL(2,\Bbb C)$.  
 
Our double group is complete since we have the following unique  
(Iwasawa) decompositions, where  
$$\align 
\ph\i: SL(2,\Bbb C) &\to SU(2).SB(2,\Bbb C),  
     \text{ where } s=\frac1{\sqrt{|z_1|^2+|z_3|^2}},\\ 
\pmatrix z_1 & z_2\\ z_3 & z_4\endpmatrix  
     &= \pmatrix sz_1 & -s\bar z_3 \\  
          sz_3 & s\bar z_1 \endpmatrix  
     \pmatrix 1/s & s(\bar z_1z_2+\bar z_3z_4) \\ 0 & s \endpmatrix, \\ 
\ps\i: SL(2,\Bbb C) &\to SB(2,\Bbb C).SU(2),  
     \text{ where }t=\frac1{\sqrt{|z_3|^2+|z_4|^2}},\\ 
\pmatrix z_1 & z_2\\ z_3 & z_4\endpmatrix  
&= \pmatrix t & t(z_1\bar z_3+z_2\bar z_4) \\  
          0 & 1/t \endpmatrix  
     \pmatrix t\bar z_4 & -t\bar z_3 \\ tz_3 & tz_4 \endpmatrix. \\ 
\endalign$$ 
Therefore, the bracket $\{\quad,\quad\}$ is globally symplectic on  
$SL(2,\Bbb C)$. This bracket is projectable on the subgroups  
$SU(2)$ and $SB(2,\Bbb C)$, and for the `left trivialization'  
$SL(2,\Bbb C)=SU(2).SB(2,\Bbb C)$ it gives us the Poisson Lie  
brackets on $SU(2)$:   
$$\alignat2 
\{\al,\bar\al\} &= -i|\nu|^2 &\qquad \{\nu,\bar\nu\} &= 0 \\ 
\{\al,\nu\} &= \tfrac12 i\al\nu &\qquad  
     \{\bar\al,\bar\nu\} &= -\tfrac12i\bar\al\bar\nu\\ 
\{\al,\bar\nu\} &= \tfrac12 i\al\bar\nu  &\qquad  
     \{\bar\al,\nu\} &=-\tfrac12i\bar\al\nu,  
\endalignat$$ 
and on $SB(2,\Bbb C)$: 
$$ 
\{\ga,t\} = \tfrac12 i\ga t,\qquad  
     \{\bar\ga,\ga\} = i\left(t^2-\frac1{t^2}\right).  
$$ 
It is possible to linearize $SB(2,\Bbb C)$ with this Lie-Poisson  
structure. The mapping  
$$\gather 
\pmatrix t & \ga\\ 0 & 1/t \endpmatrix \mapsto (\log(t),  
\operatorname{Re}\om, \operatorname{Im}\om), 
     \quad\text{ where }\\ 
\om=\sqrt{\frac{R^2-\log^2(t)}{|\ga|^2}}.\ga,\quad 
     R=\frac12\operatorname{arcosh}\left 
     (\frac{|\ga|^2 + t^2+ 1/t^2}{2}\right) 
\endgather$$ 
gives us a Poisson diffeomorphism between $(SB(2,\Bbb C),\La^{G_-})$  
and the linear Poisson structure defining the coadjoint bracket on  
$\frak s\frak u(2)$, namely  
$z\partial_x\wedge \partial_y + y\partial_z\wedge \partial_x +  
x\partial_y\wedge \partial_z$. These formulae were first obtained by  
Xu, see also \cit!{36}. 
 
Since $H^2(SL(2,\Bbb C))=0$, the symplectic structure $\om=\La_+\i$  
is exact, so there is a potential $\Th$ with $d\Th=\om$.  
Moreover, $(SL(2,\Bbb C),\La_+)$ is symplectomorphic to $T^*SU(2)$  
with the canonical symplectic structure, since  
$(G_-=SB(2,\Bbb C),\La^{G_-})$ is  
Poisson equivalent to $\frak s\frak b(2,\Bbb C)$ with its  
$\frak s\frak u(2)$-dual Poisson structure. So from the Poisson point  
of view there is no difference between $(SL(2,\Bbb C),\La_+)$ and  
$T^*SU(2)$ (they are isomorphic as symplectic groupoids), 
but the group structures differ. 
 
\subhead\nmb.{7.3}. Example \endsubhead 
On the "$ax+b$" Lie algebra $\g_+$ spanned by $X_1,X_2$ with commutator  
$[X_1,X_2]=X_2$ the cobracket given by $b'(X_1)=0$ and  
$b'(X_2)=X_1\wedge X_2$  
defines a Lie bialgebra structure. The Lie bracket on $\g_-=\g_+^*$  
is then given by $[Y_1,Y_2]=Y_2$, and the remaining commutator  
relations on $\g=\g_+\oplus\g_-$ is given by  
$[X_1,Y_1]=0$, $[X_1,Y_2]=-Y_1$, $[X_2,Y_1]=X_2$,  
$[X_2,Y_2]=-X_1+Y_1$. 
A matrix representation of $\g$ is the Lie algebra  
$\frak g\frak l(2,\Bbb R)$ via 
$$ 
X_1=\pmatrix 1 & 0 \\ 0 & 0 \endpmatrix, 
X_2=\pmatrix 0 & 1 \\ 0 & 0 \endpmatrix, 
Y_1=\pmatrix 0 & 0 \\ 0 & 1 \endpmatrix, 
Y_2=\pmatrix 0 & 0 \\ 1 & 0 \endpmatrix, 
$$ 
with the metric  
$$ 
\ga(A,B) = \tr(AJBJ),\quad\text{ where }  
     J=\pmatrix 0 & 1 \\ 1 & 0 \endpmatrix. 
$$ 
The subgroups $G_\pm$ of the Lie group $G=GL^+(2,\Bbb R)$ of matrices  
with determinant $>0$ are given by 
$$ 
G_+ = \left\{\pmatrix x & y \\ 0 & 1 \endpmatrix : x>0  \right\}, 
G_- = \left\{\pmatrix 1 & 0 \\ a & b \endpmatrix : b>0  \right\}, 
$$ 
The calculation of the affine Poisson tensor $\La_+$ on $G$ in the  
coordinates   
$$ 
\pmatrix x & y \\ a & b \endpmatrix\quad\text{ gives }\quad 
\La_+ = xy\partial_x\wedge \partial_y + ab\partial_a\wedge \partial_b  
+ xb(\partial_x\wedge \partial_b + \partial_a\wedge \partial_y). 
$$ 
It is degenerate at points with $xb=0$ and vanishes at $x=b=0$. This  
shows that $G_+G_-\ne G$. Indeed, one can easily see that $G_+G_-$  
consists of all matrices with $b\ne 0$, and $G_-G_+$ of those with  
$x\ne 0$. This should mean that the dressing vector fields are not  
all complete. Indeed,  
$\la_+(X_1)=-x(bx-ya)\partial_x$, which restricted to $G_+$ gives  
$-x^2\partial_x$, a vector field on $\Bbb R^+$ which is not complete  
since its flow is given by  
$$ 
\Fl^{-x^2\partial_x}_{x_0}(t) = \frac{x_0}{tx_0+1}. 
$$


\Refs 
 
\widestnumber\no{99} 
 
\ref 
\no \cit0{1} 
\by Alekseev, A.Yu.; Malkin, A.Z. 
\paper Symplectic structures associated to Poisson Lie groups 
\jour Comm. Math. Phys.  
\vol 162 
\yr 1994 
\pages 147--173 
\endref 
 
\ref 
\no \cit0{2} 
\by Alekseevsky, Dmitri V\.; Grabowski, Janusz;  
Marmo, Giuseppe; Michor, Peter W. 
\paper Poisson structures on the cotangent bundle of a Lie group or  
a principle bundle and their reductions 
\jour J. Math. Physics 
\vol 35 
\yr 1994 
\pages 4909--4928 
\endref 
 
\ref  
\no \cit0{3}  
\by Alekseevsky, Dmitri; Grabowski, Janusz; Marmo, Giuseppe; Michor,   
Peter W.  
\paper Completely integrable systems: a generalization 
\paperinfo to appear 
\jour  Modern Physics Letters A 
\endref  
 
\ref 
\no \cit0{4} 
\by Astrakhantsev, V.V. 
\paper A characteristic property of simple Lie algebras 
\jour Funct. Anal. Appl. 
\vol 19 
\yr 1985 
\pages 65-66 
\endref 
 
\ref 
\no \cit0{5} 
\by Astrakhantsev, V.V. 
\paper Decomposability of metrizable Lie algebras 
\jour Funct. Anal. Appl. 
\vol 12 
\yr 1978 
\pages 64-65 
\endref 
 
\ref 
\no \cit0{6} 
\by Belavin, A.A.; Drinfeld, V.G. 
\paper On the solutions of the classical Yang-Baxter equation 
\jour Funct. Anal. Appl. 
\vol 16 
\yr 1982 
\page 159 
\endref 
 
\ref 
\no \cit0{7} 
\by Belavin, A.A.; Drinfeld, V.G. 
\paper The triangle equations and simple Lie algebras 
\paperinfo Preprint of Inst of Theoretical Physics 
\yr 1982 
\endref 
 
\ref 
\no \cit0{8} 
\by Bordemann, M. 
\paper Nondegenerate invariant bilinear forms on non--associative  
algebras  
\paperinfo Pre\-print Freiburg THEP 92/3, to appear 
\jour Acta Math. Univ. Comenianae 
\endref 
 
\ref 
\no \cit0{9} 
\by Cahen, M.; Gutt, S; Rawnsley, J. 
\paper Some remarks on the classification of Poisson Lie groups 
\inbook Symplectic Geometry and Quantization 
\eds Maeda, Y.; Omori, H.; Weinstein, A. 
\bookinfo Contemporary Mathematics 179 
\publ AMS 
\publaddr Providence 
\yr 1994 
\pages 1--16 
\endref 
 
\ref 
\no \cit0{10} 
\by Drinfeld, V.I. 
\paper Hamiltonian structures on Lie groups, Lie bialgebras, and the  
geometric meaning of Yang-Baxter equations 
\jour Dokl. Akad. Nauk SSSR 
\vol 268,2 
\yr 1983 
\pages 285--287 
\endref 
 
\ref 
\no \cit0{11} 
\by Drinfeld, V.I. 
\paper Quantum groups 
\inbook Proceedings of the International Congress of Mathematicians,  
Berkeley, California, USA, 1986  
\bookinfo Vol. 1 
\publ AMS 
\publaddr  
\yr 1987 
\pages 798--820 
\endref 

\ref 
\no \cit0{12} 
\by Duflo, M; Vergne, M. 
\paper Une propriet\'e de la repr\'esentation coadjointe d'une 
alg\'ebre de Lie  
\jour C.R. Acad. Sci. Paris S\'er. A-B  
\vol 268 
\yr 1969 
\pages A583--A585 
\endref 
 
\ref 
\no \cit0{13} 
\by Grabowski, J.; Marmo, G.; Perelomov, A. 
\paper Poisson structures: towards a classification 
\jour Mod. Phys. Lett. A 
\vol 8 
\yr 1993 
\pages 1719--1733 
\endref 
 
\ref 
\no \cit0{14} 
\by Kac, V. 
\book Infinite dimensional Lie algebras 
\publ Cambridge University Press 
\yr 1990 
\endref 
 
\ref 
\no \cit0{15} 
\by Karasev, M.V. 
\paper Analogues of objects of the Lie group theory for non-linear  
Poisson brackets 
\jour Soviet Mat. Izviestia 
\vol 28 
\yr 1987 
\pages 497--527 
\endref 
 
\ref  
\no \cit0{16} 
\by Kol\'a\v r, Ivan; Slov\'ak, Jan; Michor, Peter W. 
\book Natural operations in differential geometry   
\publ Springer-Verlag 
\publaddr Berlin, Heidelberg, New~York 
\yr 1993 
\pages vi+434 
\endref 
 
\ref 
\no \cit0{17} 
\by Lecomte, P. B. A.; Roger, C.  
\paper Modules et cohomologies des bigebres de Lie  
\jour C. R. Acad. Sci. Paris  
\vol 310  
\yr 1990  
\pages  405--410  
\moreref  
\paper (Note rectificative)  
\jour C. R. Acad. Sci. Paris  
\vol 311  
\yr 1990  
\pages  893--894  
\endref 
 
\ref 
\no  \cit0{18} 
\by Lizzi, F.; Marmo, G; Sparano, G.; Vitale, P. 
\paper Dynamical aspects of Lie-Poisson structures 
\jour Mod. Phys. Lett. A 
\vol 8 
\yr 1993 
\pages 2973--2987 
\endref 
 
\ref 
\no  \cit0{19} 
\by Liu, Zhang-Ju; Qian, Min 
\paper Generalized Yang-Baxter equations, Koszul Operators and  
Poisson Lie groups 
\jour J. Diff. Geom. 
\vol 35 
\yr 1992 
\pages 399--414 
\endref 
 
\ref 
\no  \cit0{20} 
\by Lu, J-H. 
\book Multiplicative and affine Poisson structures on Lie groups 
\publ Thesis 
\publaddr Berkeley 
\yr 1990 
\endref 
 
\ref 
\no  \cit0{21} 
\by Lu, J-H.; Weinstein, A. 
\paper Poisson Lie groups, dressing transformations,  
and Bruhat decompositions 
\jour J. Diff. Geom. 
\vol 31 
\yr 1990 
\pages 501--526 
\endref 
 
\ref 
\no \cit0{22} 
\by Majid, S. 
\paper Matched pairs of Lie groups associated to solutions of the  
Yang-Baxter equations 
\jour Pac. J. Math. 
\vol 141 
\yr 1990 
\pages 311--332 
\endref 
 
\ref 
\no  \cit0{23} 
\by Marmo, G.; Simoni, A; Stern, A. 
\paper Poisson Lie group symmetries for the isotropic rotator 
\jour Int. J. Mod. Phys. A 
\vol 10 
\yr 1995 
\pages 99--114 
\endref 
 
\ref
\no  \cit0{24} 
\by Marsden, J; Ratiu, T.
\book Introduction to mechanics and symmetry
\publ Springer-Verlag
\publaddr New York, Berlin, Heidelberg
\yr 1994
\endref

\ref 
\no  \cit0{25} 
\by Medina, A.; Revoy, Ph. 
\paper Algebres de Lie et produit scalaire invariant 
\jour Ann. Sci. Ec. Norm. Super., IV. Ser.  
\vol 18 
\yr 1985 
\pages 553-561  
\endref 
 
\ref 
\no  \cit0{26} 
\by Medina, A.; Revoy, Ph. 
\paper La notion de double extension et les groupes de Lie-Poisson.  
\jour Semin. Gaston Darboux Geom. Topologie Differ.  
\vol 1987-1988 
\yr 1988 
\pages 141-171  
\endref 
 
\ref  
\no  \cit0{27} 
\by Michor, Peter W. 
\paper The cohomology of the diffeomorphism group is a Gelfand-Fuks  
cohomology  
\jour Suppl. Rendiconti del Circolo Matematico di Palermo, Serie II, 
\vol 14  
\yr 1987 
\pages 235-- 246 
\finalinfo ZB~634.57015, MR~89g:58228 
\endref 
 
\ref  
\no  \cit0{28} 
\by Michor, Peter W. 
\paper Remarks on the Schouten-Nijenhuis bracket 
\jour Suppl. Rendiconti del Circolo Mate\-matico di Palermo, Serie II, 
\vol 16  
\yr 1987 
\pages 208--215 
\finalinfo ZB~646.53013, MR~89j:58003 
\endref 
 
\ref  
\no  \cit0{29} 
\by Michor, Peter W. 
\paper Knit products of graded Lie algebras and groups 
\jour Suppl. Rendiconti Circolo Mate\-matico di Palermo, Ser. II 
\vol 22 
\yr 1989 
\pages 171--175 
\finalinfo MR~91h:17024 
\endref 
 
\ref 
\no \cit0{30} 
\by Nijenhuis, A.; Richardson, R.     
\paper Cohomology and deformations in graded Lie algebras 
\jour Bull. AMS 
\vol 72 
\yr 1966 
\pages 1--29 
\endref 
 
\ref 
\no \cit0{31} 
\by Semenov-Tian-Shansky, M.A. 
\paper What is a classical $R$-matrix 
\jour Funct. Anal. Appl. 
\vol 17, 4 
\yr 1983 
\pages 17--33 
\endref 
 
\ref 
\no \cit0{32} 
\by Semenov-Tian-Shansky, M.A. 
\paper Dressing transformations and Poisson Lie group actions 
\jour Publ RIMS 
\vol 21 
\yr 1985 
\pages 1237--1260 
\endref 
 
\ref 
\no \cit0{33} 
\by Semenov-Tian-Shansky, M.A. 
\paper Poisson-Lie groups, quantum duality principle, and the twisted  
quantum double 
\jour Theor. Math. Phys. 
\vol 93 
\yr 1992 
\pages 302--329 
\lang Russian 
\endref 
 
\ref 
\no \cit0{34}  
\by Sz\'ep, J.  
\paper On the structure of groups which can be represented as the  
product of two subgroups   
\jour Acta Sci. Math. Szeged  
\vol 12  
\yr 1950  
\pages 57--61  
\endref 
 
\ref 
\no \cit0{35} 
\by Vaisman, I. 
\book Lectures on the geometry of Poisson manifolds 
\publ Birkh\"auser 
\publaddr Boston 
\yr 1994 
\endref 
 
\ref 
\no \cit0{36} 
\by Zakrzewski, S. 
\paper Classical mechanical systems based on Poisson geometry 
\paperinfo Preprint 
\endref 
 
\ref 
\no \cit0{37} 
\by Zha, Jianguo 
\paper Fixed-point-free automorphisms of Lie algebras 
\jour Acta Math. Sin., New Ser.  
\vol 5, 1 
\yr 1989 
\pages 95-96 
\endref 
 
\endRefs
\enddocument